\documentclass[10pt,a4paper]{article}

\usepackage[a4paper, total={6in, 8in}]{geometry}

\usepackage{amsmath}
\usepackage{amsfonts}
\usepackage{amssymb}

\usepackage{graphicx}
\usepackage{amsmath}
\usepackage{amsfonts}

\usepackage{amsthm}
\usepackage{epstopdf}
\usepackage[round]{natbib}

\usepackage{lipsum} 

\usepackage{algorithm}
\usepackage{algpseudocode}

\newtheorem{theorem}{Theorem}
\newtheorem{assumption}{Assumption}
\newtheorem{lemma}{Lemma}
\newtheorem{corollary}{Corollary}

\theoremstyle{definition}
\newtheorem{example}{Example}[section]

\title{Central-limit approach to risk-aware Markov decision processes}

\author{
	Pengqian Yu\thanks{This work was supported in part by the EU FP7 project INSIGHT under grant 318225 while J. Yu and P. Yu were at IBM Research Ireland.}  \thanks{P. Yu is with the Department of Mechanical Engineering, National University of Singapore, 9 Engineering Drive 1, Singapore 117575, Singapore (e-mail: yupengqian@nus.edu.sg).} 
	\and Jia Yuan Yu\footnotemark[1] \thanks{J. Yu is with Concordia Institute of Information System Engineering, Concordia University, 1455 de Maisonneuve Blvd. West EV007.635, Montreal, Quebec H3G 1M8, Canada (e-mail: jy@osore.ca).} 
	\and Huan Xu \thanks{H. Xu is with the Department of Mechanical Engineering, National University of Singapore, 9 Engineering Drive 1, Singapore 117575, Singapore (e-mail: mpexuh@nus.edu.sg).}}

\begin{document}
	\date{}
	\maketitle

	\begin{abstract}
		Whereas classical Markov decision processes maximize the expected
		reward, we consider minimizing the risk. We propose to evaluate the
		risk associated to a given policy over a long-enough time horizon
		with the help of a central limit theorem.  The proposed approach
		works whether the transition probabilities are known or not. We also
		provide a gradient-based policy improvement algorithm that converges
		to a local optimum of the risk objective.
	\end{abstract}
	
	\section{Introduction}
	
	Markov Decision Processes (MDPs) are essential models for stochastic
	sequential decision-making problems
	\citep[e.g.,][]{puterman2014markov,bertsekas1995dynamic,sutton1998reinforcement}. Classical
	MDPs are concerned with the expected performance criteria. However, in
	many practical problems, risk-neutral objectives may not be
	appropriate \citep[cf.][Example 1]{ruszczynski2010risk}.
	
	Risk-aware decision-making is prevalent in the financial mathematics
	and optimization literature
	\citep[e.g.,][]{artzner2002coherent,rockafellar2007coherent,markowitz2000mean,von2007theory}, but
	limited to single-stage decisions. Risk-awareness has been adopted in
	multi-stage or sequential decision problems more recently. A
	chance-constrained formulation for MDPs with uncertain parameters has
	been discussed in \citep{delage2010percentile}. A criteria of
	maximizing the probability of achieving a pre-determined target
	performance has been analyzed in \citep{xu2011probabilistic}. In
	\citep{ruszczynski2010risk}, Markov risk measure is introduced and a
	corresponding risk-averse policy iteration method is
	developed. A mean-variance optimization
	problem for MDPs is addressed in \citep{mannor2013algorithmic}. A generalization of percentile optimization with objectives defined by a measure over percentiles instead of a single percentile is introduced in \citep{chen2012tractable}. In terms of
	computational convenience, actor-critic algorithms for optimizing
	variance-related risk measures in both discounted and average reward
	MDPs have been proposed in \citep{prashanth2013actor}. Mean and
	conditional value-at-risk optimization problem in MDPs is solved by
	policy gradient and actor-critic algorithms in
	\citep{chow2014algorithms}. A unified policy gradient approach is
	proposed in \citep{tamar2015policy} to seek optimal strategy in MDPs
	with the whole class of coherent risk measures.  Risk-aversion in
	multiarmed bandit problems is studied in \citep{yu2013sample}, where
	algorithms are proposed with PAC guarantees for best arm
	identification. Robust MDPs \citep[e.g.,][]{nilim2005robust,iyengar2005robust} mitigate the sensitivity of optimal policy to ambiguity in the underlying transition probabilities, and are closely related to the risk-sensitive MDPs. Such relations are uncovered by \citet{osogami2012robustness}.
	% However, we consider the approach proposed in \cite{yu2013sample} is
	% quite restrictive. First, in general, the rewards $r_i$,
	% $i = 1,\cdots,T$ generated by following a Markov chain induced by a
	% policy are usually dependent. Secondly, as $T\rightarrow\infty$, the
	% method proposed in \cite{yu2013sample} is intractable. Specifically,
	% the density estimation of $r_1+\cdots+r_\tau$ by convolution becomes
	% computationally hard.
	
	In this work, we study the value-at-risk
	in finite-horizon MDPs.  Computing this risk associated with a sequence of decision
	using the defintion and first principles is intractable.  However, we
	show that this computation is made tractable by using a central limit
	theorem for Markov chains
	\citep{jones2004markov,kontoyiannis2003spectral,glynn1996liapounov}.
	For a long-enough horizon $T$ and a fixed policy $\pi$, we
	are thus able to evaluate the risk associated with following policy
	$\pi$ over $T$ time steps.
	% This risk can be expressed as a function of the distribution of
	% cumulative reward.
	Specifically, our first contributions are policy evaluation algorithms
	whether the transition probability matrix induced by $\pi$ is known or
	not. Under mild conditions, we provide high-probability error bounds
	for the evaluation.  
	
	For a fixed risk measure $\rho$, and a space of policies that is
	parametrized by $\theta$, our second contribution is a policy improvement algorithm
	that converges in finite iterations, under certain assumptions, to a
	locally optimal policy.  This approach updates the parameter $\theta$
	in the direction of gradient of the risk measure.  Compared to the
	previous work, our proposed method does not explicitly approximate the
	value function. Therefore, it does not have approximation error due to
	the selection of value function approximator.
	
	Even though we deal mainly with the static value-at-risk
	risk measure (defined by the cumulative reward), our results can be easily extended to 
	deal with the conditional- or average-value-at-risk by using its definition \cite{Schied04}.

	\subsection{Distinction from related works}
	
	It is important to note that our risk measure is
	very different from the dynamic risk measure analyzed in \citep{ruszczynski2010risk},
	which has the following recursive structure: for $T$-length horizon MDP and
	policy $\pi$, the dynamic risk measure $\rho_T$ is defined
	as
	$$r(x_0^\pi)+\gamma\rho\left(r(x_1^\pi) +\gamma\rho\Big(r(x_2^\pi)+ \gamma\rho\big(\ldots +r(x_T^\pi)\big)\Big)\right),$$
	where $r$ is a state dependent reward function, $\gamma<1$ is a constant and $\rho$ is a
	Markov risk measure, i.e., the evaluation of each static coherent risk
	measure $\rho$ is not allowed to depend on the whole past.
	In contrast, our risk
	measure emphasizes the statistical properties of the reward accumulated
	over multiple time steps: it is similar to that used in
	\citep{yu2013sample} in the context of bandit problems.

	From the computational and algorithmic perspective, we evaluate the policy and obtain the risk value directly. Different from the work \citep{prashanth2013actor,chow2014algorithms,tamar2015policy}, which indirectly evaluated the risk value by using function approximation, our proposed method has an explicit policy evaluation step and does not require value function approximation. Therefore, it does not have approximation error due to the selection of value function approximator and the richness of the features. Moreover, for risk-aware MDPs that are known to be NP-hard \citep[e.g.,][]{delage2010percentile,xu2011probabilistic}, our method can find approximated solutions to those problems. 
	
	From the conceptual perspective, our work is the first attempt to consider general risk measures with a central-limit approach. Compared with \citep{tamar2015policy}, our approach is not limited to so-called ``coherent'' risk measures (both static and dynamic), whose defining properties are not satisfied by the value-at-risk. The methods proposed in \citep{tamar2015policy} thus cannot be applied. Moreover, the time horizon in our setup is allowed to be infinite. Other previous work \citep{delage2010percentile,xu2011probabilistic,mannor2013algorithmic} restrict themselves to specific risk criteria: \citet{delage2010percentile} considered a set of percentile criteria; \citet{xu2011probabilistic} seek to find the policy that maximizes the probability of achieving a pre-determined target performance, and \citet{mannor2013algorithmic} discussed a mean-variance optimization problem.
	
	This paper is organized as follows: In Section 2 we provide background and problem statement. For a fixed policy, we then evaluate its risk value in Section 3. In specific, we discuss the policy evaluation algorithms in Section 3.1 and 3.2 when the transition kernel is either known or unknown. The main result is Theorem \ref{thm_N}, which bounds the error of true and estimated asymptotic variances. In Section 4, we provide policy gradient algorithm to improve the fixed policy. A conclusion is offered in Section 5. 
	
	\section{Problem formulation}
	
	MDP is a tuple $\langle\mathcal{X},\mathcal{A},P,r,T\rangle$, where
	$\mathcal{X}$ is a finite set of states, $\mathcal{A}$ is a finite set
	of actions, $P$ is the transition kernel, $T$ is the (possibly
	infinite) moderate time horizon, and
	$r(x,a):\mathcal{X}\times\mathcal{A}\rightarrow\mathbb{R}$ is the
	bounded and measurable reward function. A function $F:\mathcal{X}\rightarrow\mathbb{R}$ is called lattice if there are $h>0$ and $0\leq d<h$, such that $\frac{F(x)-d}{h}$ is an integer for $x\in\mathcal{X}$. The minimal $h$ for such condition holds is called the span of $F$. If the function $F$ can be written as a sum $F=F_0+F_l$, where $F_l$ is lattice with span $h$ and $F_0$ has zero asymptotic variance (which will be defined later), then $F$ is called almost-lattice. Otherwise, $F$ is called strongly nonlattice. We assume the functional $r$
	is strongly nonlattice. A policy $\pi$ is the rule according to which
	actions are selected at each state. We parameterize the stochastic
	policies by
	$\{\pi(\cdot|x;\theta),x\in\mathcal{X},\theta\in\Theta\subseteq\mathbb{R}^{k_1}\}$,
	and denote the set of policies by $\Pi$. Since in this setting a
	policy $\pi$ is represented by its $k_1$-dimensional parameter vector
	$\theta$, policy dependent functions can be written as a function of
	$\theta$ in place of $\pi$. So, we use $\pi$ and $\theta$
	interchangeably in the paper.
	
	Given a policy $\pi\in\Pi$, a Markov reward process
	$X^\pi=\{x_t^\pi:x\in\mathcal{X},t\in T\}$ is induced by $\pi$.  We
	write the transition kernel of the induced Markov process and {\em
		dependent} rewards as $P^\pi$ and $r(x^\pi_t)$ for all $t\in T$,
	respectively. We denote $(P^\pi)^n$ as the $n$-step Markov transition
	kernel corresponding to $P^\pi$. We are interested in the long-term
	behavior of the sum
	\begin{equation*}
		S_T^\pi=\sum_{t=0}^{T-1}r(x_t^\pi),\quad T\geq1.
	\end{equation*}
	If the Markov chain is positive recurrent with invariant probability
	measure $\xi^\pi$, we define the mean of reward function $r$ as
	\begin{equation}\label{eq_average reward}
		\varphi^\pi(r)=\lim_{T\rightarrow\infty}\frac{1}{T}\mathbb{E}\left[\sum_{t=0}^{T-1}r(x_t^\pi)\right],
	\end{equation}
	where $x_0^\pi$ is the initial condition and $\mathbb{E}$ is the
	expectation over the Markov process $X^\pi$. Equivalently, the mean
	can be expressed as
	\begin{equation}\label{mean}
		\varphi^\pi(r)=\sum_{x}\xi^\pi(x)r(x).
	\end{equation}
	Often we can quantify the rate of convergence of \eqref{eq_average
		reward} by showing that the limit
	\begin{equation*}
		\hat{r}^\pi(x)=\lim_{T\rightarrow\infty}\mathbb{E}\left[\sum_{t=0}^{T-1}r(x_t^\pi)-T\varphi^\pi(r)\right]
	\end{equation*}
	exists where, in fact, the function $\hat{r}^\pi(x)$ solves the
	Poisson equation
	\begin{equation}\label{eq_poisson}
		P^\pi\hat{r}^\pi=\hat{r}^\pi-r+\varphi^\pi(r)\mathbf{1}.
	\end{equation}	
	Here $P^\pi$ operates on functions
	$\hat{r}^\pi:\mathcal{X}\rightarrow\mathbb{R}$ via
	$P^\pi\hat{r}^\pi(x)=\sum_{y}P(x,y)\hat{r}^\pi(y)$, $r$ is the vector
	form of $r(x)$ and $\mathbf{1}$ is the vector with all elements being
	$1$.
	
	Let $\mathcal{Y}$ denote the set of bounded random variables, and let
	$\rho:\mathcal{Y}\rightarrow\mathbb{R}$ denote a risk measure. Our
	risk-aware objective is
	\begin{equation}\label{objective}
		\min_{\pi\in\Pi}\quad\rho\left(S_T^\pi\right).
	\end{equation}
	We discuss specific types of risk measures in the coming sections. All
	proofs are deferred in the supplementary material due to space
	constraints.
	
	\section{Policy evaluation}
	
	In this section, we develop methods to evaluate the risk value
	$\rho(S^\pi_T)$ for a fixed policy $\pi$. In specific, we consider the
	Markov chain $X^\pi$ induced by policy $\pi$, which is assumed to be
	geometrically ergodic. We then apply limit theorems for geometrically
	ergodic Markov chain in \citep{kontoyiannis2003spectral} to obtain the
	distribution function of cumulative reward when transition kernel is
	either known or unknown. Finally, we provide algorithms for computing
	some examples of risk measures.
	
	Similar to \citet{kontoyiannis2003spectral}, we make the following two
	standard assumptions for the Markov process $X^\pi$.
	
	\begin{assumption}\label{asu_ergodicity}
		The Markov process $X^\pi$ is geometrically ergodic with Lyapunov
		function $V$ (i.e., $\psi-$irreducible, aperiodic and satisfying
		Condition (V4) of \citet{kontoyiannis2003spectral}) where
		$V:\mathcal{X}\rightarrow[1,\infty)$ satisfies
		$\varphi^\pi(V^2)<\infty$.
	\end{assumption}
	The detailed definition of geometric ergodicity can be found in \citep[Section 2]{kontoyiannis2003spectral}. 
	\begin{assumption}\label{asu_reward}
		The (measurable) reward function $r:\mathcal{X}\rightarrow[-1,1]$
		has nontrivial asymptotic variance
		$\sigma^2\triangleq\lim_{t\rightarrow\infty}var_{x}(S^\pi_t/\sqrt{t})>0$.
	\end{assumption}
	
	Under Assumption \ref{asu_ergodicity} and \ref{asu_reward},
	\citet{kontoyiannis2003spectral} showed that the following result
	holds.
	
	\begin{theorem}[Theorem 5.1, \citet{kontoyiannis2003spectral}]\label{thm_cdf}
		Suppose that $X^\pi$ and the strongly nonlattice functional $r$
		satisfy Assumption \ref{asu_ergodicity} and \ref{asu_reward}, and
		let $G^\pi_T(y)$ denote the distribution function of the normalized
		partial sums $\frac{S^\pi_T-T\varphi^\pi(r)}{\sigma\sqrt{T}}$:
		\begin{equation*}
			G^\pi_T(y)\triangleq\mathbb{P}\left\{\frac{S^\pi_T-T\varphi^\pi(r)}{\sigma\sqrt{T}}\leq y \right\}.
		\end{equation*}
		Then, for all $x^\pi_0\in\mathcal{X}$ and as $T\rightarrow\infty$,
		\begin{equation*}
			G^\pi_T(y)=g(y)+\frac{\gamma(y)}{\sigma\sqrt{T}}\left[\frac{\varrho}{6\sigma^2}(1-y^2)-\hat{r}^\pi(x^\pi_0)\right]+o(T^{-\frac{1}{2}}),
		\end{equation*}
		uniformly in $y\in\mathbb{R}$, where $\gamma(y)$ denotes the
		standard normal density and $g(y)$ is the corresponding distribution
		function. Here $\varrho$ is a constant related to the third moment
		of $S_T^\pi/\sqrt{T}$. The definitions of $\varrho$ and term
		$o(T^{-1/2})$ can be found in the supplementary material.
	\end{theorem}
	
	Unlike \citep{heyman2003stochastic}, which considered finite Markov chains and obtained the moments of the cumulative reward, we further impose geometric ergodicity assumption on the Markov process $X^\pi$. By doing this, we are able to utilize the explicit form of distribution function given by \citep[Theorem 5.1]{kontoyiannis2003spectral} with the time error term $o(T^{-1/2})$. In addition, the geometric ergodicity property enables us to bound the approximation error in Theorem \ref{thm_N}. 
	
	Similar to \citet{glynn1996liapounov}, we define the {\em fundamental
		kernel} $Z^\pi$ associated with fixed $\pi$ as
	\begin{equation*}
		Z^\pi=(H^\pi)^{-1},
	\end{equation*}
	where the kernel $\Xi^\pi$ is defined as
	$\Xi^\pi(x,\cdot)=\xi^\pi(\cdot)$, and
	$H^\pi\triangleq I-P^\pi-\Xi^\pi$. If the inverse of $H^\pi$ is well
	defined, \citet{glynn1996liapounov} stated that the solution to
	Poisson equation \eqref{eq_poisson} has the form
	\begin{equation}\label{sol_poisson}
		\hat{r}^\pi=Z^\pi(r-\varphi^\pi(r)\mathbf{1}).
	\end{equation}
	
	When $\sigma^2<\infty$ and $\varphi^\pi(\hat{r}^\pi)<\infty$,
	\citet[Equation 17.44]{meyn2012markov} showed that the asymptotic
	variance $\sigma^2$ can be obtained by
	\begin{equation}\label{asymptotic variance}
		\begin{aligned}
			\sigma^2&=\lim_{n\rightarrow\infty}\frac{\mathbb{E}[(S^\pi_n-n\varphi^\pi(r))^2]}{n}=\sum_{x\in\mathcal{X}}[\hat{r}^\pi(x)^2-(P^\pi\hat{r}^\pi(x))^2]\xi^\pi(x).
		\end{aligned}
	\end{equation}
	We remark that the asymptotic variance can also be calculated without
	the knowledge of $\hat{r}^\pi$ by Theorem 17.5.3 of
	\citet{meyn2012markov}:
	\begin{equation*}
		\begin{aligned}
			\sigma^2=\sum_{x^\pi_0\in\mathcal{X}}{\overline{r}_0}^2\xi^\pi(x^\pi_0)+2\sum_{k=1}^\infty\sum_{x^\pi_0\in\mathcal{X}}\overline{r}_0\xi^\pi(x^\pi_0)\times\sum_{x^\pi_k\in\mathcal{X}}\overline{r}_kP^\pi(x,x^\pi_k|x=x^\pi_0)^k,
		\end{aligned}
	\end{equation*}
	where $\overline{r}_k=r(x^\pi_k)-\varphi^\pi(r)$ and
	$P^\pi(x,y|x=x_0)^n$ denotes the $n$-step transition probability to
	state $y$ conditioned on the initial state $x_0$. It is clear that
	solving \eqref{asymptotic variance} requires the knowledge of
	transition kernel $P^\pi$. In the following subsections, we study the
	policy evaluation with and without knowing $P^\pi$, respectively.
	
	\subsection{Transition probability kernel $P^\pi$ is known}
	
	In this subsection, we propose an algorithm to obtain the distribution
	function $G_T^\pi$ in Theorem \ref{thm_cdf}, using which we can
	evaluate the risk value.
	
	Since transition kernel $P^\pi$ is known, we first calculate the
	stationary distribution $\xi^\pi$. Note that the stationary
	distribution satisfies
	\begin{equation*}
		P^\pi\xi^\pi=\xi^\pi.
	\end{equation*}
	Thus, the stationary distribution is the eigenvector of matrix $P^\pi$
	with eigenvalue $1$ (recall that $1$ is always an eigenvalues for
	matrix $P^\pi$), which is unique up to constant multiples under
	Assumption \ref{asu_ergodicity}.
	
	After we compute the stationary distribution $\xi^\pi$, it is straightforward to calculate the mean \eqref{mean}, solve Poisson equation
	\eqref{eq_poisson} and get the asymptotic variance \eqref{asymptotic
		variance}. Therefore, we are able to obtain the distribution
	function $G_T^\pi(y)$ defined in Theorem \ref{thm_cdf} at this
	stage. Finally, the risk measures that can be expressed as a function
	of $G_T^\pi(y)$ can be evaluated. We provide two examples below.
	\begin{example}[$T$-step value-at-risk]
		Let $\lambda$ be given and fixed. The $T$-step value-at-risk is
		defined by
		\begin{equation*}
			\rho_\lambda(S^\pi_T)=VaR_{\lambda}(S_T^\pi)=-q(\lambda),
		\end{equation*}
		where $q$ is the right-continuous quantile
		function\footnote{Formally,
			$q(\lambda)=\inf\{y\in\mathbb{R}:F(y)>\lambda\}$, where $F$ is the
			distribution function of $x$.} of $S_T^\pi$.
	\end{example}
	When $T$ is moderate, $o(T^{-1/2})$ vanishes, the distribution
	function $G^\pi_T(y)$ defined in Theorem \ref{thm_cdf} has the form
	\begin{equation}\label{estimated VaR_P known}
		G^\pi_T(y)\approx g(y)+\frac{\gamma(y)}{\sigma\sqrt{T}}\left[\frac{\varrho}{6\sigma^2}(1-y^2)-\hat{r}^\pi(x^\pi_0)\right].
	\end{equation}
	This allows us to compute the $T$-step value-at-risk. The procedure is
	summarized in Algorithm \ref{algo_policy evaluation_P known}. The next
	example considers the mean-variance risk measure.
	\begin{example}[Mean-variance risk]
		Given $\lambda$, the mean-variance risk measure is given by
		\begin{equation*}
			\rho_\lambda(S^\pi_T)=-\mu(S^\pi_T)+\lambda\sigma^2(S^\pi_T/\sqrt{T}),
		\end{equation*}
		where $\mu(S^\pi_T)$ is the mean of $S^\pi_T$, and
		$\sigma^2(S^\pi_T/\sqrt{T})$ is the variance of $S^\pi_T/\sqrt{T}$.
	\end{example}
	The corollary below gives formula to obtain the mean-variance risk. An
	algorithm can be obtained similarly for the mean-variance risk simply
	by replacing line $7$ in Algorithm \ref{algo_policy evaluation_P
		known} by the corresponding risk measure. The same procedure can
	easily be generalized to other risk measures.
	\begin{corollary}\label{col_1}
		Under Assumption \ref{asu_ergodicity} and \ref{asu_reward}, given
		$\lambda$, the mean-variance risk measure can be computed by
		\begin{equation*}
			\rho_\lambda(S^\pi_T)=-\varphi(r)+\lambda\sigma^2.
		\end{equation*}
		where $\varphi(r)$ and $\sigma^2$ are defined in \eqref{mean} and
		\eqref{asymptotic variance}, respectively.
	\end{corollary}
	\begin{algorithm}[H]
		\caption{Policy evaluation for $T$-step value-at-risk ($P^\pi$ is
			known)}
		\label{algo_policy evaluation_P known}
		\begin{algorithmic}[1]
			\State \textbf{Input:} Transition probability matrix $P^\pi$
			induced by policy $\pi$, initial state $x^\pi_0$, decision horizon $T\in\mathbb{N}_+$,
			reward function $r$, and the coefficient $\lambda$ of
			value-at-risk.  \State Obtain stationary distribution $\xi^\pi$ by
			solving the equation $P^\pi\xi^\pi=\xi^\pi$.  \State Substitute
			$\xi^\pi$ in \eqref{mean} for mean
			$\varphi^\pi(r)=\sum_{x}\xi^\pi(x)r(x)$.  \State Solve the Poisson
			equation \eqref{eq_poisson} for $\hat{r}^\pi$ by
			\eqref{sol_poisson}.  \State Compute asymptotic variance in
			\eqref{asymptotic variance}:
			$\sigma^2=\sum_{x\in\mathcal{X}}[\hat{r}^\pi(x)^2-(P^\pi\hat{r}^\pi(x))^2]\xi^\pi(x)$.
			\State Calculate the constant $\varrho$.  \State \textbf{Output:}
			The $T$-step value-at-risk, which can be computed by
			\begin{equation*}
				\widetilde{VaR}_{\lambda}(S_T^\pi)=-\inf\left\{y\in\mathbb{R}:G^\pi_T(h(y))>\lambda\right\},
			\end{equation*}
			where $G^\pi_T(y)$ is defined in \eqref{estimated VaR_P known}, and $h(y)=\frac{y-T\varphi^\pi(r)}{\sigma\sqrt{T}}$.
		\end{algorithmic}
	\end{algorithm}
	
	\subsection{Transition probability kernel $P^\pi$ needs to be
		estimated}
	
	In practice, knowing transition kernel $P^\pi$ is generally hard. In
	this subsection, we present a policy evaluation algorithm when
	transition kernel needs to be estimated, and provide the corresponding
	approximation error.
	
	Intuitively, if the estimated transition kernel is accurate enough, the gap between solutions to estimated and true Poisson equations is small (as shown in Lemma \ref{lem_solution difference}), which implies the difference between estimated and true asymptotic variances is also small (as shown in Theorem \ref{thm_N}).
	
	A way of estimating transition kernel is to use the empirical
	distribution. Denoting the required number of samples by $n_1$ and the
	estimated transition kernel by $P_{n_1}^\pi$, we make the following
	assumption.
	\begin{assumption}\label{asu_estimated probability matrix}
		There exist $\epsilon_1(n_1)>0$ and $\delta_1(n_1)>0$ such that
		$\epsilon_1(n_1)\rightarrow0$ and $\delta_1(n_1)\rightarrow0$ as
		$n_1\rightarrow\infty$. Moreover, we have the following bound for
		the error\footnote{We use $\epsilon_1(n_1)$ and $\delta_1(n_1)$ to
			show the dependence of the error on $n_1$.}:
		\begin{equation*}
			\mathbb{P}\left(||P^\pi_{n_1}-P^\pi||_{\infty}\leq\epsilon_1(n_1)\right)\geq 1-\delta_1(n_1).
		\end{equation*} 		
	\end{assumption}
	Here the supreme norm $||\cdot||_{\infty}$ of a matrix $A$ is defined
	as the maximum absolute row sum of the matrix:
	\begin{equation*}
		||A||_\infty\triangleq\max_{1\leq i\leq m}\sum_{j=1}^{n}|a_{ij}|.
	\end{equation*}
	We remark Assumption \ref{asu_estimated probability matrix} can be easily satisfied if a large collection of observed data is provided. This could be achievable if one is allowed to simulate the system arbitrarily many times.
	
	Next, we perturb the true transition kernel $P^\pi$, define the
	resulting perturbed transition kernel as $\tilde{P}^\pi$ and denote
	the stationary distribution associated with $\tilde{P}^\pi$ as
	$\tilde{\xi}^\pi$. Given $\tilde{P}^\pi$, we let
	$\xi^\pi_n(\cdot)=\tilde{P}^\pi(x,\cdot)^n$ be an estimator of
	$\tilde{\xi}^\pi(\cdot)$. Under Assumption \ref{asu_ergodicity}, we
	make the following assumption which establishes the convergence of
	$||\xi^\pi_n-\tilde{\xi}^\pi||$, where $||\cdot||$ is the total
	variation norm defined as following: If $\mu$ is a signed measure on
	$\mathcal{B}(\mathcal{X})$, then the total variation norm is defined
	as
	\begin{equation*}
		||\mu||=\sup_{|f|\leq1}|\mu(f)|=\sup_{A\in\mathcal{B}(\mathcal{X})}\mu(A)-\inf_{A\in\mathcal{B}(\mathcal{X})}\mu(A).
	\end{equation*}
	\begin{assumption}\label{asu_convergence rate}
		Given $\epsilon>0$ and perturbed transition kernel
		$\tilde{P}^\pi\in P$. If
		$||\tilde{P}^\pi-P^\pi||_\infty\leq\epsilon,$ then for every
		$n\in\mathbb{N}_+$, the geometrically ergodic Markov chain has a
		form of convergence
		\begin{equation*}
			||\tilde{P}^\pi(x,\cdot)^n-\tilde{\xi}^\pi(\cdot)||\leq M(x)\tau^{n},
		\end{equation*}
		where $\tau<1$, and $M(x)$ is a nonnegative function.
	\end{assumption}
	This assumption essentially guarantees the geometric ergodicity
	property remains if the original Markov chain is slightly perturbed.
	
	Assumption \ref{asu_estimated probability matrix} and
	\ref{asu_convergence rate} immediately imply that
	\begin{equation}\label{probabilistic convergence rate}
		\mathbb{P}\left(||P_{n_1}^\pi(x,\cdot)^{n}-\tilde{\xi}^\pi(\cdot)||\leq
		M(x)\tau^{n}\right)\geq1-\delta_1(n_1),
	\end{equation}
	where $\tilde{\xi}^\pi$ is the stationary distribution associated with
	$P^\pi_{n_1}$. 
	Let $\tau_1(P^\pi)$ be the ergodicity coefficient of transition
	kernel $P^\pi$ \citep{seneta1988perturbation}, which is defined by
	\begin{equation}\label{def_ergodicity coefficient}
		\tau_1(P^\pi)=\sup_{\substack{||\nu||_1=1\\\nu^\top \mathbf{e}=0}}||\nu^\top P^\pi||_1.
	\end{equation}
	If the ergodicity coefficient satisfies $0\leq\tau_1(P^\pi)<1$,
	\citet[Equation 3.4]{cho2001comparison}
	shows
	$$||\tilde{\xi}^\pi-\xi^\pi||_1\leq\frac{||P^\pi_{n_1}-P^\pi||_{\infty}}{1-\tau_1(P^\pi)}.$$
	Therefore, under Assumption \ref{asu_estimated probability matrix}, we
	have
	\begin{equation}\label{invariant probabilities convergence}
		\mathbb{P}\left(||\tilde{\xi}^\pi-\xi^\pi||_1\leq\frac{\epsilon_1(n_1)}{1-\tau_1(P^\pi)}\right)\geq1-\delta_1(n_1).
	\end{equation}	
	The ergodicity coefficient $\tau_1(P_{n_1}^\pi)$ for the estimated transition kernel $P_{n_1}^\pi$ is defined in a same manner.
	
	Matrix $H^\pi$ plays an important role in defining solution to Poisson
	equation. We impose a requirement on the smallest singular value
	$\sigma_{min}(H^\pi)$ of $H^\pi$.
	\begin{assumption}\label{asu_poisson eq}
		There exists a constant $c>\widetilde{\varpi}(x)$ such that the
		smallest singular value of matrix $H^\pi$ satisfies
		$\sigma_{min}(H^\pi)\geq c$. Here\footnote{The parameter $\widetilde{\varpi}$ is
			dependent on the initial state $x$ through nonnegative function
			$M(x)$. In the following, we omit $x$ for convenience and emphasize
			$x$ when such dependence is required to show.}
		\begin{equation}\label{def_tilde w}
			\begin{aligned}
				\widetilde{\varpi}(x)&=\sqrt{|\mathcal{X}|}\left(\epsilon_1(n_1)+\frac{\epsilon_1(n_1)}{1-\tau_1(P^\pi_{n_1})-\epsilon_1(n_1)}+2M(x)\tau^{n_2}\right).
			\end{aligned}
		\end{equation}
	\end{assumption}
	Similar to \eqref{eq_poisson}, we define the estimated Poisson
	equation as
	\begin{equation}\label{eq_estimated poisson}
		P^\pi_{n_1}\hat{r}_{n_1,n_2}^\pi=\hat{r}_{n_1,n_2}^\pi-r+\varphi_{n_2}^\pi(r)\mathbf{1},
	\end{equation}
	where $\varphi_{n_2}^\pi(r)$ is the estimated mean
	\begin{equation}\label{estimated mean}
		\varphi_{n_2}^\pi(r)=\sum_{x}\xi^\pi_{n_2}(x)r(x),
	\end{equation}
	with $\xi^\pi_{n_2}$ being the estimation of stationary distribution
	$\tilde{\xi}^\pi$
	\begin{equation}\label{estimated stationary distribution}
		\xi^\pi_{n_2}(\cdot)=P_{n_1}^\pi(x,\cdot)^{n_2}.
	\end{equation} 	
	We then define the {\em perturbed fundamental kernel} as
	$$Z_{n_1,n_2}^\pi=(H_{n_1,n_2}^\pi)^{-1},$$
	where $H_{n_1,n_2}^\pi=I-P_{n_1}^\pi-\Xi_{n_2}^\pi$ with
	$\Xi_{n_2}^\pi(x,\cdot)=\xi_{n_2}^\pi(\cdot)$. The following lemma
	shows $Z_{n_1,n_2}^\pi$ is well defined under Assumption
	\ref{asu_poisson eq}.
	\begin{lemma}\label{lem_Z12}
		Under Assumption \ref{asu_poisson eq}, $Z_{n_1,n_2}^\pi$ is well
		defined.
	\end{lemma}
	The solution to estimated Poisson equation \eqref{eq_estimated
		poisson} can be expressed as
	\begin{equation}\label{sol_estimated poisson}
		\hat{r}_{n_1,n_2}^\pi=Z_{n_1,n_2}^\pi(r-\varphi_{n_2}^\pi(r)\mathbf{1}),
	\end{equation}
	
	We define the asymptotic variance associated with $\hat{r}_{n_1,n_2}$
	as $\sigma_{n_1,n_2}^2$, which can be obtained by the equation below:
	\begin{equation}\label{estimated asymptotic variance}
		\sigma_{n_1,n_2}^2=
		\sum_{x\in\mathcal{X}}[\hat{r}_{n_1,n_2}^\pi(x)^2-(P_{n_1}^\pi\hat{r}_{n_1,n_2}^\pi(x))^2]\xi_{n_2}^\pi(x).
	\end{equation}	
	The theorem below provides the condition for $n_2$ such that the
	difference between estimated and true asymptotic variances $\sigma^2$
	and $\sigma_{n_1,n_2}^2$ defined in \eqref{asymptotic variance} and
	\eqref{estimated asymptotic variance} is small. We provide an outline
	of the proof below and a complete version can be found in the
	supplementary material.
	
	\begin{theorem}\label{thm_N}	
		Given policy $\pi$, $\epsilon\geq0$, $\lambda\in(0,1)$ and let
		$\sigma_{n_1,n_2}^2$ and $\sigma^2$ be the estimated and true
		asymptotic variances, respectively. Under Assumption
		\ref{asu_ergodicity}, \ref{asu_reward}, \ref{asu_estimated
			probability matrix}, \ref{asu_convergence rate} and
		\ref{asu_poisson eq}, there exists $n_2$ such that the following two
		conditions (C1) and (C2) are satisfied
		\begin{equation*}
			\begin{aligned}
				\text{(C1)}\qquad
				&\frac{|\mathcal{X}|(|\mathcal{X}|+1)^3(\widetilde{\varpi}-\sqrt{|\mathcal{X}|}\epsilon_1(n_1))}{(c-\widetilde{\varpi})^2}\leq\lambda\epsilon\\
				\text{(C2)}\qquad
				&\frac{\alpha_2(x)|\mathcal{X}|}{c}\left(\frac{(|\mathcal{X}|+1)\widetilde{\varpi}}{c-\widetilde{\varpi}}+\sqrt{|\mathcal{X}|}\left(\widetilde{\varpi}-\epsilon_1(n_1)\sqrt{|\mathcal{X}|}\right)\right)+\max_{x\in\mathcal{X}}\alpha_1(x)\alpha_2(x)|\mathcal{X}|\leq(1-\lambda)\epsilon,
			\end{aligned}
		\end{equation*}
		where
		\begin{equation*}
			\begin{aligned}
				\alpha_1(x)=&\frac{\epsilon_1(n_1)\sqrt{|\mathcal{X}|}(|\mathcal{X}|+1)}{c-\widetilde{\varpi}}+\frac{|\mathcal{X}|^{\frac{3}{2}}}{c}\left(\frac{\widetilde{\varpi}(|\mathcal{X}|+1)}{c-\widetilde{\varpi}}+\sqrt{|\mathcal{X}|}\left(\widetilde{\varpi}-\epsilon_1(n_1)\sqrt{|\mathcal{X}|}\right)\right),\\
				\alpha_2(x)=&\dfrac{|\mathcal{X}|(|\mathcal{X}|+1)(2c-\widetilde{\varpi})}{c(c-\widetilde{\varpi})}.\end{aligned}
		\end{equation*}
		Moreover, we
		have
		$$\mathbb{P}(|\sigma_{n_1,n_2}^2-\sigma^2|\leq\epsilon)\geq1-38\delta_1(n_1).$$
	\end{theorem}
	\begin{proof}
		By equations \eqref{asymptotic variance} and \eqref{estimated
			asymptotic variance}, we have
		\begin{equation*}
			\begin{aligned}	
				&|\sigma_{n_1,n_2}^2-\sigma^2|\\
				=&\bigg|\sum_{x\in\mathcal{X}}[{\hat{r}_{n_1,n_2}^\pi}(x)^2-(P_{n_1}^\pi\hat{r}_{n_1,n_2}^\pi(x))^2]\xi^\pi_{n_2}(x)-\sum_{x\in\mathcal{X}}[{\hat{r}^\pi}(x)^2-(P^\pi\hat{r}^\pi(x))^2]\xi^\pi(x)\bigg|\\
				\stackrel{(1)}{=}&\bigg|\sum_{x\in\mathcal{X}}[{\hat{r}_{n_1,n_2}^\pi}(x)^2-(P_{n_1}^\pi\hat{r}_{n_1,n_2}^\pi(x))^2]\xi^\pi_{n_2}(x)-\sum_{x\in\mathcal{X}}[{\hat{r}_{n_1,n_2}^\pi}(x)^2-(P_{n_1}^\pi\hat{r}_{n_1,n_2}^\pi(x))^2]\xi^\pi(x)\\
				&+\sum_{x\in\mathcal{X}}[{\hat{r}_{n_1,n_2}^\pi}(x)^2-(P_{n_1}^\pi\hat{r}_{n_1,n_2}^\pi(x))^2]\xi^\pi(x)-\sum_{x\in\mathcal{X}}[{\hat{r}^\pi}(x)^2-(P^\pi\hat{r}^\pi(x))^2]\xi^\pi(x)\bigg|\\
				&\stackrel{(2)}{\leq}a_n+b_n,
			\end{aligned}
		\end{equation*}
		where
		\begin{equation*}
			\begin{aligned}
				a_n=&\bigg|\sum_{x\in\mathcal{X}}[{\hat{r}_{n_1,n_2}^\pi}(x)^2-(P_{n_1}^\pi\hat{r}_{n_1,n_2}^\pi(x))^2](\xi^\pi_{n_2}(x)-\xi^\pi(x))\bigg|,\\
				b_n=&\bigg|\sum_{x\in\mathcal{X}}[{\hat{r}_{n_1,n_2}^\pi}(x)^2-(P_{n_1}^\pi\hat{r}_{n_1,n_2}^\pi(x))^2]\xi^\pi(x)-\sum_{x\in\mathcal{X}}[{\hat{r}^\pi}(x)^2-(P^\pi\hat{r}^\pi(x))^2]\xi^\pi(x)\bigg|.
			\end{aligned}
		\end{equation*}
		Equality $(1)$ holds due to subtracting and then adding a same term
		$\sum_{x\in\mathcal{X}}[{\hat{r}_{n_1,n_2}^\pi}(x)^2-(P_{n_1}^\pi\hat{r}_{n_1,n_2}^\pi(x))^2]\xi^\pi(x)$. Inequality
		$(2)$ holds by triangle inequality. In order to bound the terms
		$a_n$ and $b_n$, we need following lemmas.
		
		First, we derive a bound for $|\tau_1(P_{n_1}^\pi)-\tau_1(P^\pi)|$ by the
		lemma below.
		\begin{lemma}\label{lem_ergodicity coefficient}
			Under Assumption \ref{asu_estimated probability matrix}, the
			ergodicity coefficient for the estimated transition kernel
			$P_{n_1}^\pi$
			satisfies
			$$\mathbb{P}\left(|\tau_1(P_{n_1}^\pi)-\tau_1(P^\pi)|\leq\epsilon_1(n_1)\right)\geq1-\delta_1(n_1).$$
		\end{lemma}
		
		Denoting $||a||_2$ as the Euclidean norm of vector $a$, we present a
		lemma showing that the solution \eqref{sol_poisson} of Poisson
		equation \eqref{eq_poisson} is bounded.
		\begin{lemma}\label{lem_hat(r)}
			Under Assumption \ref{asu_reward} and \ref{asu_poisson eq}, the
			solution to Poisson equation $\hat{r}^\pi$ is bounded
			$||\hat{r}^\pi||_2\leq\frac{\sqrt{|\mathcal{X}|}(|\mathcal{X}|+1)}{c}$.
		\end{lemma}
		
		Analogously to Lemma \ref{lem_hat(r)}, we have the following result.
		\begin{lemma}\label{lem_hat(r)_n}
			Under Assumption \ref{asu_ergodicity}, \ref{asu_reward},
			\ref{asu_estimated probability matrix}, \ref{asu_convergence rate}
			and \ref{asu_poisson eq}, the solution \eqref{sol_estimated
				poisson} to estimated Poisson equation \eqref{eq_estimated
				poisson} is bounded such that
			$$\mathbb{P}\left(||\hat{r}_{n_1,n_2}^\pi||_{2}\leq\frac{\sqrt{|\mathcal{X}|}(|\mathcal{X}|+1)}{c-\widetilde{\varpi}}\right)\geq1-4\delta_1(n_1).$$
		\end{lemma}
		The lemma below shows that solutions to the true and estimated Poisson equations are close to each
		other with high probability.
		
		\begin{lemma}\label{lem_solution difference}
			Under Assumption \ref{asu_ergodicity}, \ref{asu_reward},
			\ref{asu_estimated probability matrix}, \ref{asu_convergence
				rate} and \ref{asu_poisson eq}, the difference between
			solutions to true and estimated Poisson equation is bounded
			such that
			\begin{equation*}
				\mathbb{P}\left(||\hat{r}^\pi-\hat{r}_{n_1,n_2}^\pi||_2\leq\tilde{\varsigma}\right)\geq1-7\delta_1(n_1),
			\end{equation*}
			where
			$\tilde{\varsigma}=\frac{\sqrt{|\mathcal{X}|}(|\mathcal{X}|+1)\widetilde{\varpi}}{c(c-\widetilde{\varpi})}+\frac{|\mathcal{X}|}{c}\left(\widetilde{\varpi}-\epsilon_1(n_1)\sqrt{|\mathcal{X}|}\right)$.
		\end{lemma}

		Recall that for two random variables $a_n$ and $b_n$, and
		$\epsilon\in\mathbb{R}$ with $\lambda\in(0,1)$, by union bound
		and monotonicity property of probability, we have
		\begin{equation*}
			\mathbb{P}(a_n+b_n\leq \epsilon)\geq\mathbb{P}(a_n\leq\lambda \epsilon)+\mathbb{P}(b_n\leq(1-\lambda) \epsilon).
		\end{equation*}
		Finally, by bounding terms $a_n$ and $b_n$, we prove the
		theorem.
	\end{proof}
	Now we are able to compute the $T$-step value-at-risk, which is
	summarized as a procedure in the Algorithm \ref{algo_policy
		evaluation_P unknown}. Note that given the asymptotic variance
	$\sigma_{n_1,n_2}$ and solution to estimated Poisson equation
	$\hat{r}_{n_1,n_2}$, when $T$ is moderate, $o(T^{-1/2})$
	vanishes, $G^\pi_T(y)$ defined in Theorem \ref{thm_cdf} becomes
	\begin{equation}\label{estimated VaR_P unknown}
		G^\pi_T(y)\approx g(y)+\frac{\gamma(y)}{\sigma_{n_1,n_2}\sqrt{T}}\left[\frac{\tilde{\varrho}}{6\sigma_{n_1,n_2}^2}(1-y^2)-\hat{r}_{n_1,n_2}^\pi(x^\pi_0)\right],
	\end{equation}
	where $\tilde{\varrho}$ is defined in a same manner as $\varrho$
	by replacing $P^\pi$ and $\xi^\pi$ by $P^\pi_{n_1}$ and
	$\xi^\pi_{n_2}$, respectively.
	
	In the following, we analysis the mean-variance risk defined in
	Example 3.2. Corollary \ref{col_2} gives formula to obtain the
	mean-variance risk and provides high-probability bound for the
	risk measure. An algorithm can be developed for the
	mean-variance risk simply by replacing line $9$ in Algorithm
	\ref{algo_policy evaluation_P unknown} by the corresponding risk
	measure.
	\begin{corollary}\label{col_2}
		Under Assumption \ref{asu_ergodicity}, \ref{asu_reward},
		\ref{asu_estimated probability matrix}, \ref{asu_convergence
			rate} and \ref{asu_poisson eq}, given $\lambda$, the
		estimated mean-variance risk measure can be computed by
		\begin{equation*}
			\tilde{\rho}_\lambda(S^\pi_T)=-\varphi_{n_2}^\pi(r)+\lambda\sigma_{n_1,n_2}^2.
		\end{equation*}
		where $\varphi_{n_2}^\pi(r)$ and $\sigma_{n_1,n_2}^2$ are
		defined in \eqref{estimated mean} and \eqref{estimated
			asymptotic variance}, respectively. Moreover, we have
		\begin{equation*}
			\mathbb{P}\left(|\tilde{\rho}_\lambda(S^\pi_T)-\rho_\lambda(S^\pi_T)|\leq\kappa\right)\geq1-41\delta_1(n_1),
		\end{equation*}
		where $\rho_\lambda(S^\pi_T)$ is the true mean-variance risk
		value computed in Corollary \ref{col_1} and
		$\kappa=2M(x)\tau^{n_2}+\frac{\epsilon_1(n_1)}{1-\tau_1(P_{n_1}^\pi)-\epsilon_1(n_1)}+\lambda\epsilon$.
	\end{corollary}
	A similar argument can be developed for $T$-step value-at-risk defined in
	Example 3.1. When $\varrho=1$, $\hat{r}=0.5$, $T=100$, $\varphi(r)=0.1$, we plot the distribution functions $G_T$ for $\sigma=1$ and $1.1$ below.
	\begin{figure}[ht]
		\centering
		\includegraphics[width=0.5\textwidth]{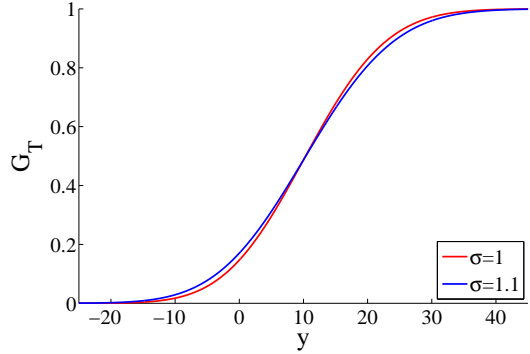}
		\caption{Distribution functions $G_T$ for $\sigma=1$ and $1.1$.}\label{fig1}
	\end{figure}
	Figure \ref{fig1} shows the distribution functions are close, which implies the corresponding $T$-step value-at-risks are also close to each other (When $\lambda=0.3$, the $T$-step value-at-risks are $5.134$ and $4.632$, respectively).   
	
	%	--- version 2 ---
	\begin{algorithm}[H]
		\caption{Policy evaluation for $T$-step value-at-risk ($P^\pi$
			is unknown)}
		\label{algo_policy evaluation_P unknown}
		\begin{algorithmic}[1]
			\State \textbf{Input:} Policy $\pi$, initial state $x^\pi_0$, reward function $r$,
			decision horizon $T\in\mathbb{N}_+$, $\lambda\in(0,1)$ and
			estimation parameters $\epsilon_1(n_1)$ and $\delta_1(n_1)$.
			\State Obtain the estimated transition kernel $P^\pi_{n_1}$
			satisfying Assumption \ref{asu_estimated probability matrix}
			with parameter $\epsilon_1(n_1)$ and $\delta_1(n_1)$.
			\State Compute the ergodicity coefficient
			$\tau_1(P^\pi_{n_1})$ for the estimated transition kernel
			$P^\pi_{n_1}$.  \State Choose $n_2$ with
			$\tau_1(P^\pi_{n_1})$ obtained from step $3$ such that
			conditions in Theorem \ref{thm_N} hold.  \State Compute the
			estimated stationary distribution $\xi_{n_2}^\pi$ defined by
			\eqref{estimated stationary distribution}.  \State Calculate
			the solution of estimated Poisson equation
			$\hat{r}_{n_1,n_2}^\pi$ by \eqref{sol_estimated poisson}.
			\State Obtain the asymptotic variance $\sigma_{n_1,n_2}^2$
			defined in \eqref{estimated asymptotic variance}.  \State
			Calculate the constant $\tilde{\varrho}$.  \State
			\textbf{Output:} The $T$-step value-at-risk, which can be
			computed by
			\begin{equation*}
				\widetilde{VaR}_{\lambda}(S_T^\pi)=-\inf\left\{y\in\mathbb{R}:G^\pi_T(h(y))>\lambda\right\},
			\end{equation*}
			where $G^\pi_T(y)$ is defined in \eqref{estimated VaR_P unknown}, and
			$h(y)=\frac{y-T\varphi_{n_2}^\pi(r)}{\sigma_{n_1,n_2}\sqrt{T}}$
			with $\varphi_{n_2}^\pi(r)$ defined in \eqref{estimated
				mean}.
			% \EndProcedure
		\end{algorithmic}
	\end{algorithm}
	\section{Policy Improvement}
	
	In this section, we propose policy gradient methods to improve
	policy $\pi$ and solve the optimization problem
	\eqref{objective}. Specifically, we compute the gradient of the
	performance (i.e., the risk) w.r.t.\ the policy parameters from
	the evaluated policies, and improve $\pi$ by adjusting its
	parameters in the direction of the gradient.
	
	In the following, we make a standard assumption in
	gradient-based MDPs literature
	\citep[e.g.,][]{sutton1999policy,konda1999actor,bhatnagar2009natural,prashanth2013actor}.
	
	\begin{assumption}\label{asu_policy parameter}
		For any state-action pair $(x,a)$, $\pi(a|x;\theta)$ is continuously
		differentiable in the parameter $\theta$.
	\end{assumption}
	
	In a recent work, \citet[Section 4]{prashanth2013actor} presented
	actor-critic algorithms for optimizing the risk-sensitive measure that
	are based on two simultaneous perturbation methods: simultaneous
	perturbation stochastic approximation (SPSA)
	\citep{spall1992multivariate} and smoothed functional (SF)
	\citep{bhatnagar2012stochastic}. We use similar arguments here and
	further apply those methods to estimate the gradient of the risk
	measure with respect to the policy parameter $\theta$, i.e.,
	$\nabla_{\theta}\rho(S_T^\theta)$.
	
	The idea in these methods is to estimate the gradient
	$\nabla_{\theta}\rho(S_T^\theta)$ using two simulated trajectories of
	the system corresponding to policies with parameters $\theta$ and
	$\theta^+=\theta+\beta\triangle$. Here $\beta>0$ is a positive
	constant and $\triangle$ is a perturbation random variable, i.e., a
	$k_1$-vector of independent Rademacher (for SPSA) and Gaussian
	$\mathcal{N}(0,1)$ (for SF) random variables.
	
	SPSA-based estimate for $\nabla_{\theta}\rho(S_T^\theta)$ is given by
	\begin{equation}\label{g_1}
		\partial_{\theta^{(i)}}\rho(S_T^\theta)\approx\frac{\rho(S_T^{\theta+\beta\triangle})-\rho(S_T^\theta)}{\beta\triangle^{(i)}},\quad i=1,\dots,k_1
	\end{equation}
	where $\triangle$ is a vector of independent Rademacher random variables
	and $\triangle^{(i)}$ is its $i$-th entry. The advantage of this
	estimator is that it perturbs all directions at the same time (the
	numerator is identical in all $k_1$ components). So, the number of
	function measurements needed for this estimator is always two,
	independent of the dimension $k_1$.
	
	SF-based estimate for $\nabla_{\theta}\rho(S_T^\theta)$ is given by
	\begin{equation}\label{g_2}
		\partial_{\theta^{(i)}}\rho(S_T^\theta)\approx\frac{\triangle^{(i)}}{\beta}\left(\rho(S_T^{\theta+\beta\triangle})-\rho(S_T^\theta)\right),i=1,\dots,k_1
	\end{equation}
	where $\triangle$ is a vector of independent Gaussian
	$\mathcal{N}(0,1)$ random variables.
	
	The step size of the gradient descent $\{\lambda(t)\}$ satisfies the
	following condition.
	\begin{assumption}\label{asu_step size}
		The step size schedule $\{\lambda(t)\}$ satisfies
		\begin{equation*}
			\lambda(t)>0,\quad\sum_{t}\lambda(t)=\infty,\quad\text{and}\quad\sum_{t}\lambda(t)^2<\infty.
		\end{equation*}
	\end{assumption}
	At each time step $t$, we update the policy parameter $\theta$ as
	follows: for $i=1,\cdots,k_1$,
	\begin{equation}\label{update}
		\theta_{t+1}^{(i)}=\Gamma_i\left(\theta^{(i)}_t-\lambda(t)\partial_{\theta^{(i)}}\rho(S_T^\theta)\right).
	\end{equation}
	Note that $\triangle_t^{(i)}$'s are independent Rademacher and
	Gaussian $\mathcal{N}(0,1)$ random variables in SPSA and SF updates,
	respectively. $\Gamma$ is an operator that projects a vector
	$\theta\in\mathbb{R}^{k_1}$ to the closest point in a compact and
	convex set $C\subset\mathbb{R}^{k_1}$. The projection operator is
	necessary to ensure the convergence of the algorithms. Specifically,
	we consider the differential equation
	\begin{equation}\label{differential equation}
		\dot{\theta}_t=\check{\Gamma}(\nabla_{\theta}\rho(S_T^{\theta_t})),
	\end{equation}
	where $\check{\Gamma}$ is defined as follows: for any bounded
	continuous function $f(\cdot)$,
	$$\check{\Gamma}(f(\theta_t))=\lim_{\tau\rightarrow0}\frac{\Gamma(\theta_t+\tau f(\theta_t))-\theta_t}{\tau}.$$
	The projection operator $\check{\Gamma}(\cdot)$ ensures that the
	evolution of $\theta$ via the differential equation
	\eqref{differential equation} stays within the bounded set
	$C\in\mathbb{R}^{k_1}$. Let
	$\mathcal{Z}=\{\theta\in C:
	\check{\Gamma}(\nabla_{\theta}\rho(S_T^{\theta_t}))=0\}$
	denote the set of asymptotically stable equilibrium points of the
	differential equation \eqref{differential equation} and
	$\mathcal{Z}^\epsilon=\{\theta\in C:
	||\theta-\theta_0||<\epsilon,\theta_0\in\mathcal{Z}\}$
	denote the set of points in the $\epsilon$-neighborhood of
	$\mathcal{Z}$.
	
	The following is our main result, which bounds the error of our policy improvement algorithm.
	
	\begin{theorem}\label{thm_improvement}
		Under Assumption \ref{asu_ergodicity}, \ref{asu_policy parameter}
		and \ref{asu_step size}, for any given $\epsilon>0$, there exists
		$\beta_0>0$ such that for all $\beta\in(0,\beta_0)$,
		$\theta_t\rightarrow\theta^*\in\mathcal{Z}^\epsilon$ almost surely.
	\end{theorem}
	
	The above theorem guarantees the convergence of the SPSA and SF
	updates to a local optimum of the risk objective function
	$\rho(S_T^\theta)$. Finally, we provide the policy improvement method
	in Algorithm \ref{algo_policy improvement}.
	
	\begin{algorithm}[H]
		\caption{Policy improvement}
		\label{algo_policy improvement}
		\begin{algorithmic}[1]
			\State Initialize policy $\pi$ parameterized by parameter
			$\theta_1$ and stopping criteria $\epsilon$. Let $t=1$ and $\theta_0=\theta_1+2\epsilon\mathbf{1}$. 
			\While {$||\theta_{t}-\theta_{t-1}||_\infty>\epsilon$} 
			\State Evaluate the risk
			value $\rho(S_T^{\theta_t})$ and
			$\rho(S_T^{\theta_t+\beta\triangle_t})$ associated with policy
			$\theta_t$ and $\theta_t+\beta\triangle_t$ by Algorithm
			\ref{algo_policy evaluation_P known} or \ref{algo_policy
				evaluation_P unknown}. Here $\triangle_t$ are independent
			Rademacher and Gaussian $\mathcal{N}(0,1)$ random variables in
			SPSA and SF updates, respectively.  
			\State For $i=1,\cdots,k_1$, update the policy parameter
			$$\theta_{t+1}^{(i)}=\Gamma_i\left(\theta^{(i)}_t-\lambda(t)\partial_{\theta^{(i)}}\rho(S_T^\theta)\right),$$ where $\lambda(t)$ is the step size and
			$\partial_{\theta^{(i)}}\rho(S_T^\theta)$ is defined in
			\eqref{g_1} for SPSA or \eqref{g_2} for SF.
			\State $t=t+1$.
			\EndWhile
			\State \textbf{Output:} Policy parameter $\theta$.
		\end{algorithmic}
	\end{algorithm}
	
	\section{Conclusion}
	
	In the context of risk-aware Markov decision processes,
	we use a central limit theorem efficiently evaluate the risk
	associated to a given policy. We also provide a policy
	improvement algorithm, which works on a parametrized policy-space in a
	gradient-descent way. Under mild conditions, it is guaranteed that the policy
	evaluation is approximately correct and that the policy improvement
	converges to a local optimum of the risk measure. We believe that many
	important problems that are usually addressed using standard MDP
	models should be revisited using our proposed
	approach to risk management (e.g., machine replacement,
	inventory management, portfolio optimization, etc.).
	
	% \subsubsection*{Acknowledgements}

	\appendix
	\section{Appendix}
	\subsection{Defintion of $\varrho$ in Theorem \ref{thm_cdf}}
	\begin{flalign*}
		&\varrho=\varrho_1+\varrho_3+\varrho_3,&
	\end{flalign*}
	where
	\begin{flalign*}
		&\varrho_1=\sum_{x^\pi_0\in\mathcal{X}}r(x^\pi_0)^3\xi^\pi(x^\pi_0),&
	\end{flalign*}
	\begin{flalign*}
		&\varrho_2=
		3\sum_{\substack{i=-\infty\\i\neq0}}^{\infty}\sum_{\substack{x^\pi_0\in\mathcal{X}}}r(x^\pi_0)^2\xi^\pi(x^\pi_0)\times\sum_{\substack{x^\pi_i\in\mathcal{X}}}r(x^\pi_i)P^\pi(x,x^\pi_i|x=x^\pi_0)^i,&
	\end{flalign*}
	and
	\begin{flalign*}
		&\varrho_3=
		6\sum_{i,j=1}^{\infty}\sum_{\substack{x^\pi_0\in\mathcal{X}}}r(x^\pi_0)\xi^\pi(x^\pi_0)\times\sum_{\substack{x^\pi_i\in\mathcal{X}}}r(x^\pi_i)P^\pi(x,x^\pi_i|x=x^\pi_0)^i,\times\sum_{\substack{x^\pi_j\in\mathcal{X}}}r(x^\pi_{i+j})P^\pi(x,x^\pi_j|x=x^\pi_i)^j.&
	\end{flalign*}
	Here $P^\pi(x,y|x=x_0)^n$ denotes the $n$-step transition probability
	to state $y$ conditioned on the initial state $x_0$.
	\subsection{Definition of $o(T^{-1/2})$ in Theorem \ref{thm_cdf}}
	\citet{kontoyiannis2003spectral} showed that the term $o(T^{-1/2})$ can
	be represented as
	\begin{equation*}
		o(T^{-1/2})=\frac{1}{\pi}\int_{-A\sqrt{T}}^{A\sqrt{T}}\bigg|M_T\left(\frac{i\omega}{\sigma\sqrt{T}}\right)-\phi_T(\omega)\bigg|\frac{d\omega}{|\omega|}+\frac{\nu}{\sqrt{T}}
	\end{equation*}
	where
	$$M_T(\alpha)\triangleq m_T(\alpha)\exp(-\alpha\mathbb{E}[S_T^\pi]),$$
	with $m_T(\alpha)=\mathbb{E}[\exp(\alpha S_T^\pi)]$, $T\geq1$ and
	$\alpha\in\mathbb{C}$. The characteristic function $\phi_T(\cdot)$ is
	defined as
	$$\phi_T(\omega)\triangleq\exp\left(\frac{-\omega^2}{2}\right)\left(1+\frac{\varrho(i\omega)^3}{6\sigma^3\sqrt{T}}\right),\quad\omega\in\mathbb{R}.$$	 
	$A$ is chosen large enough so that
	$A>24(\nu\pi)^{-1}|\Psi_T^{'}\left(y-\frac{\mathbb{E}[S_T^\pi]}{\sigma\sqrt{T}}\right)|$
	for all $y\in\mathbb{R}$, $T\geq1$ and arbitrary $\nu>0$. Here,
	$$\Psi_T(y)=g\left(y\right)+\frac{\varrho}{6\sigma^3\sqrt{T}}(1-y^2)\gamma\left(y\right),\quad y\in\mathbb{R}.$$
	\subsection{Proof of Lemma \ref{lem_Z12}}
	\begin{proof}
		Let $\widetilde{\Xi}^\pi(x,\cdot)=\tilde{\xi}^\pi(\cdot)$ and
		\begin{equation*}
			\begin{aligned}
				&\triangle_1=P^\pi-P_{n_1}^\pi,\\
				&\triangle_2=\Xi^\pi-\widetilde{\Xi}^\pi,\\
				&\triangle_3=\widetilde{\Xi}^\pi-\Xi_{n_2}^\pi,\\
			\end{aligned}
		\end{equation*}
		and $||A||_{2}$ be the spectral norm of matrix $A$, which is the
		largest singular value of matrix $A$:
		$$||A||_{2}=\sigma_{max}(A).$$
		
		Due to Assumption \ref{asu_estimated probability matrix},
		\eqref{probabilistic convergence rate} and \eqref{invariant
			probabilities convergence}, we have
		\begin{equation*}
			\begin{aligned}
				&\mathbb{P}\left(||\triangle_1||_{2}\leq\epsilon_1(n_1)\sqrt{|\mathcal{X}|}\right)\geq1-\delta_1(n_1),\\
				&\mathbb{P}\left(||\triangle_2||_{2}\leq\frac{\epsilon_1(n_1)\sqrt{|\mathcal{X}|}}{1-\tau_1(P^\pi)}\right)\geq1-\delta_1(n_1),\\
				&\mathbb{P}\left(||\triangle_3||_{2}\leq
				2\sqrt{|\mathcal{X}|}M(x)\tau^{n_2}\right)\geq1-\delta_1(n_1).
			\end{aligned}
		\end{equation*}
		Define $\triangle=\triangle_1+\triangle_2+\triangle_3$, we
		have
		\begin{equation}\label{triangle difference}
			\mathbb{P}\left(||\triangle||_{2}\leq\varpi(x)\right)\geq1-3\delta_1(n_1),
		\end{equation}
		where\footnote{The parameter $\varpi$ is dependent on the
			initial state $x$ through nonnegative function $M(x)$. In
			the following, we omit $x$ for convenience and emphasize $x$
			when such dependence is required to show. Same argument goes
			for $\widetilde{\varpi}$.}
		$$\varpi(x)=\sqrt{|\mathcal{X}|}\left(\epsilon_1(n_1)+\frac{\epsilon_1(n_1)}{1-\tau_1(P^\pi)}+2M(x)\tau^{n_2}\right).$$	
		We perform the singular value decomposition for the matrix
		$H^\pi$ and have
		$$H^\pi=U\Sigma V^*=U\Sigma^{1/2}I\Sigma^{1/2}V^*,$$ where $U$
		and $V^*$ are unitary matrix with dimension
		$|\mathcal{X}|\times|\mathcal{X}|$, and $\Sigma$ is a
		$|\mathcal{X}|\times|\mathcal{X}|$ diagonal matrix with
		non-negative real numbers on the diagonal. Next, by the
		definition of $Z_{n_1,n_2}^\pi$ and $\triangle$, we obtain
		\begin{equation*}
			\begin{aligned}
				Z_{n_1,n_2}^\pi=&(H^\pi+\triangle)^{-1}\\
				=&(U\Sigma^{1/2}I\Sigma^{1/2}V^*+U\Sigma^{1/2}\Sigma^{-1/2}U^{-1}\triangle (V^*)^{-1}\Sigma^{-1/2}\Sigma^{1/2}V^*)^{-1}\\
				=&(\Sigma^{1/2}V^*)^{-1}(I+\widetilde{\triangle})^{-1}(U\Sigma^{1/2})^{-1},
			\end{aligned}
		\end{equation*}	
		where
		$\widetilde{\triangle}\triangleq\Sigma^{-1/2}U^{-1}\triangle
		(V^*)^{-1}\Sigma^{-1/2}$. In addition, we have
		\begin{equation*}
			\begin{aligned}
				||\widetilde{\triangle}||_{2}&=||\Sigma^{-1}||_{2}||\triangle||_{2}.
			\end{aligned}
		\end{equation*}
		Since
		$||\Sigma^{-1}||_{2}=||(H^\pi)^{-1}||_{2}=\frac{1}{\sigma_{min}(H^\pi)}$,
		under Assumption \ref{asu_poisson eq}, we obtain
		\begin{equation}\label{tilde triangle}
			\begin{aligned}
				\mathbb{P}\left(||\widetilde{\triangle}||_{2}\leq\frac{\varpi}{c}\right)\geq1-3\delta_1(n_1).
			\end{aligned}
		\end{equation}
		Since $\varpi<\widetilde{\varpi}<c$,
		$(I+\widetilde{\triangle})^{-1}$ is well defined, which
		implies $Z_{n_1,n_2}^\pi$ is well defined.
		
	\end{proof}
	\subsection{Proof of Lemma \ref{lem_ergodicity coefficient}}
	\begin{proof}
		From the definition of $\tau_1(P^\pi)$ and
		$\tau_1(P_{n_1}^\pi)$, we have
		\begin{equation*}
			\begin{aligned}
				&|\tau_1(P_{n_1}^\pi)-\tau_1(P^\pi)|\\
				\leq&\bigg|\sup_{\substack{||\nu||_1=1\\\nu^\top \mathbf{e}=0}}||\nu^\top P_{n_1}^\pi||_1-\sup_{\substack{||\nu||_1=1\\\nu^\top \mathbf{e}=0}}||\nu^\top P^\pi||_1\bigg|\\
				\stackrel{(1)}{\leq}&\sup_{\substack{||\nu||_1=1\\\nu^\top \mathbf{e}=0}}\bigg|||\nu^\top P_{n_1}^\pi||_1-||\nu^\top P^\pi||_1\bigg|\\
				\stackrel{(2)}{\leq}&\sup_{\substack{||\nu||_1=1\\\nu^\top \mathbf{e}=0}}||\nu^\top (P_{n_1}^\pi-P^\pi)||_1\\
				\stackrel{(3)}{\leq}&\sup_{\substack{||\nu||_1=1\\\nu^\top \mathbf{e}=0}}||\nu||_1 ||P_{n_1}^\pi-P^\pi||_\infty\\
				=&||P_{n_1}^\pi-P^\pi||_\infty
			\end{aligned}
		\end{equation*}
		Here recall that if $f,g:A\rightarrow\mathbb{R}$ are bounded
		functions, then
		$\bigg|\sup_Af-\sup_Ag\bigg|\leq\sup_A|f-g|$. For
		$||\nu||_1=1$ and $\nu^\top\mathbf{e}=0$, since
		$||\nu^\top P_{n_1}^\pi||_1\leq1$ and
		$||\nu^\top P^\pi||_1\leq1$, thus inequality $(1)$
		holds. Inequality $(2)$ holds due to reverse triangle
		property, while H\"{o}lder's inequality is applied to $(3)$.
		
		From Assumption \ref{asu_estimated probability matrix}, we get
		$$\mathbb{P}\left(|\tau_1(P_{n_1}^\pi)-\tau_1(P^\pi)|\leq\epsilon_1(n_1)\right)\geq1-\delta_1(n_1).$$
	\end{proof}
	\subsection{Proof of Lemma \ref{lem_hat(r)}}
	\begin{proof}
		Recall the solution can be bounded above by
		\begin{equation*}
			\begin{aligned}
				||\hat{r}^\pi||_2&=||Z^\pi(r-\varphi^\pi(r)\mathbf{1})||_2\\
				&\leq||Z^\pi||_{sp}||r-\varphi^\pi(r)\mathbf{1}||_2\\
				&\stackrel{(1)}{\leq}\frac{\sqrt{|\mathcal{X}|}(|\mathcal{X}|+1)}{c}.
			\end{aligned}
		\end{equation*}
		Here $(1)$ is due to
		$||r-\varphi^\pi(r)\mathbf{1}||_2\leq\sqrt{|\mathcal{X}|}(|\mathcal{X}|+1)$
		and Assumption \ref{asu_poisson eq}. Since for
		$x\in\mathcal{X}$,
		\begin{equation*}
			\begin{aligned}
				|r(x)-\varphi^\pi(r)|&=\bigg|r(x)-\sum_{x'\in\mathcal{X}}\xi^\pi(x')r(x')\bigg|\\
				&\leq|r(x)|+\sum_{x'\in\mathcal{X}}\xi^\pi(x')|r(x')|\\
				&\leq|\mathcal{X}|+1,
			\end{aligned}
		\end{equation*}
		we have
		$||r-\varphi^\pi(r)\mathbf{1}||_{\infty}\leq|\mathcal{X}|+1,$
		which implies
		$||r-\varphi^\pi(r)\mathbf{1}||_2\leq\sqrt{|\mathcal{X}|}(|\mathcal{X}|+1)$.
		It is straightforward to show $||P^\pi\hat{r}^\pi||_2$ is also
		bounded.
		\begin{equation*}
			\begin{aligned}
				||P^\pi\hat{r}^\pi||_2&\leq||P^\pi||_2||\hat{r}^\pi||_2\\
				&\leq\sqrt{|\mathcal{X}|}||P^\pi||_\infty||\hat{r}^\pi||_2\\
				&\leq\frac{|\mathcal{X}|(|\mathcal{X}|+1)}{c}.
			\end{aligned}
		\end{equation*}
	\end{proof}
	\subsection{Proof of Lemma \ref{lem_hat(r)_n}}
	\begin{proof}
		Recall the solution to estimated Poisson equation can be
		expressed as
		\begin{equation*}
			\begin{aligned}
				||\hat{r}_{n_1,n_2}^\pi||_2&=||Z_{n_1,n_2}^\pi(r-\varphi_{n_2}^\pi(r)\mathbf{1})||_2\\
				&\leq||Z_{n_1,n_2}^\pi||_{2}||r-\varphi_{n_2}^\pi(r)\mathbf{1}||_2.\\
			\end{aligned}
		\end{equation*}
		First, note that
		$||H^\pi-H^\pi_{n_1,n_2}||_{2}=||\triangle||_{2}.$ From
		\eqref{triangle difference} in the proof of Lemma
		\ref{lem_Z12}, we have
		\begin{equation*}
			\mathbb{P}\left(||H^\pi-H_{n_1,n_2}^\pi||_{2}\leq\varpi\right)\geq1-3\delta_1(n_1).
		\end{equation*}
		By Weyl's inequality, we have
		\begin{equation*}
			\mathbb{P}\left(|\sigma_{min}(H^\pi)-\sigma_{min}(H_{n_1,n_2}^\pi)|\leq\varpi\right)\geq1-3\delta_1(n_1).
		\end{equation*}
		which implies
		$$\mathbb{P}\left(||Z_{n_1,n_2}^\pi||_{2}\leq\frac{1}{c-\varpi}\right)\geq1-3\delta_1(n_1).$$
		Moreover, for $x\in\mathcal{X}$,
		\begin{equation*}
			\begin{aligned}
				|r(x)-\varphi_{n_2}^\pi(r)|&=\bigg|r(x)-\sum_{x'\in\mathcal{X}}\xi_{n_2}^\pi(x')r(x')\bigg|\\
				&\leq|r(x)|+\sum_{x'\in\mathcal{X}}\xi_{n_2}^\pi(x')|r(x')|\\
				&\leq|\mathcal{X}|+1,
			\end{aligned}
		\end{equation*}
		we have
		$||r-\varphi_{n_2}^\pi(r)\mathbf{1}||_{\infty}\leq|\mathcal{X}|+1.$
		Thus
		$||r-\varphi_{n_2}^\pi(r)\mathbf{1}||_2\leq\sqrt{|\mathcal{X}|}(|\mathcal{X}|+1)$.
		Therefore, we obtain
		$$\mathbb{P}\left(||\hat{r}_{n_1,n_2}^\pi||_{2}\leq\frac{\sqrt{|\mathcal{X}|}(|\mathcal{X}|+1)}{c-\varpi}\right)\geq1-3\delta_1(n_1).$$
		By Lemma \ref{lem_ergodicity coefficient}, we have
		\begin{equation}\label{tilde_w}
			\mathbb{P}(\varpi\leq\widetilde{\varpi})\geq1-\delta_1(n_1).
		\end{equation}
		Applying \eqref{tilde_w}, we have
		$$\mathbb{P}\left(||\hat{r}_{n_1,n_2}^\pi||_{2}\leq\frac{\sqrt{|\mathcal{X}|}(|\mathcal{X}|+1)}{c-\widetilde{\varpi}}\right)\geq1-4\delta_1(n_1).$$
		It is straightforward to show
		$||P^\pi_{n_1}\hat{r}_{n_1,n_2}^\pi||_2$ is bounded because
		\begin{equation*}
			\begin{aligned}
				||P^\pi_{n_1}\hat{r}_{n_1,n_2}^\pi||_2\leq&||P^\pi_{n_1}||_2||\hat{r}_{n_1,n_2}^\pi||_2\\
				\leq&\sqrt{|\mathcal{X}|}||P^\pi_{n_1}||_\infty||\hat{r}_{n_1,n_2}^\pi||_2.\\
			\end{aligned}
		\end{equation*}
		Thus, we have
		$$\mathbb{P}\left(||P^\pi_{n_1}\hat{r}_{n_1,n_2}^\pi||_2\leq\frac{|\mathcal{X}|(|\mathcal{X}|+1)}{c-\widetilde{\varpi}}\right)\geq1-4\delta_1(n_1).$$
	\end{proof}
	\subsection{Proof of Lemma \ref{lem_solution difference}}
	\begin{proof}
		From the proof of Lemma \ref{lem_Z12}, recall that
		$Z^\pi_{n_1,n_2}$ can be expressed as
		\begin{equation*}
			\begin{aligned}
				Z_{n_1,n_2}^\pi=(\Sigma^{1/2}V^*)^{-1}(I+\widetilde{\triangle})^{-1}(U\Sigma^{1/2})^{-1}.
			\end{aligned}
		\end{equation*}	
		Since $||\widetilde{\triangle}||_{2}<1$ under Assumption
		\ref{asu_poisson eq}, we have
		$$(I+\widetilde{\triangle})^{-1}=\sum_{n=0}^{\infty}(-1)^n\widetilde{\triangle}^n.$$
		Furthermore, we have
		\begin{equation*}
			\begin{aligned}
				&||Z_{n_1,n_2}^\pi-Z^\pi||_{2}\\
				=&\bigg|\bigg|(\Sigma^{1/2}V^*)^{-1}\left(\sum_{n=1}^{\infty}(-1)^n\widetilde{\triangle}^n\right)(U\Sigma^{1/2})^{-1}\bigg|\bigg|_{2}\\
				\leq&||(\Sigma^{1/2}V^*)^{-1}||_{2}||(U\Sigma^{1/2})^{-1}||_{2}\sum_{n=1}^{\infty}||\widetilde{\triangle}||^n_{2}\\
				\leq&\frac{||\widetilde{\triangle}||_{2}}{1-||\widetilde{\triangle}||_{2}}||\Sigma^{-1}||_{2}.
			\end{aligned}
		\end{equation*}	
		Under Assumption \ref{asu_poisson eq} and \eqref{tilde
			triangle} in the proof of Lemma \ref{lem_Z12}, we have
		\begin{equation*}
			\mathbb{P}\left(||Z^\pi-Z_{n_1,n_2}^\pi||_{2}\leq\frac{\widetilde{\varpi}}{c(c-\widetilde{\varpi})}\right)\geq1-4\delta_1(n_1).
		\end{equation*}	
		
		The difference of solutions to true and estimated Poisson
		equation can be expressed as
		\begin{equation*}
			\begin{aligned}
				&||\hat{r}^\pi-\hat{r}_{n_1,n_2}^\pi||_2\\
				=&||Z^\pi(r-\varphi^\pi(r)\mathbf{1})-|Z_{n_1,n_2}^\pi(r-\varphi_{n_2}^\pi(r)\mathbf{1})||_2\\
				=&||Z^\pi(r-\varphi^\pi(r)\mathbf{1})-Z^\pi(r-\varphi_{n_2}^\pi(r)\mathbf{1})+Z^\pi(r-\varphi_{n_2}^\pi(r)\mathbf{1})-Z_{n_1,n_2}^\pi(r-\varphi_{n_2}^\pi(r)\mathbf{1})||_2\\
				=&||(Z^\pi-Z_{n_1,n_2}^\pi)(r-\varphi_{n_2}^\pi(r)\mathbf{1})+Z^\pi(\varphi_{n_2}^\pi(r)-\varphi^\pi(r))\mathbf{1}||_2\\
				\leq&||(Z^\pi-Z_{n_1,n_2}^\pi)(r-\varphi_{n_2}^\pi(r)\mathbf{1})||_2+||Z^\pi(\varphi_{n_2}^\pi(r)-\varphi^\pi(r))\mathbf{1}||_2\\
				\leq&||Z^\pi-Z_{n_1,n_2}^\pi||_{2}||r-\varphi_{n_2}^\pi(r)\mathbf{1}||_2+||Z^\pi||_{2}||(\varphi_{n_2}^\pi(r)-\varphi^\pi(r))\mathbf{1}||_2\\
				\leq&||Z^\pi-Z_{n_1,n_2}^\pi||_{2}||r-\varphi_{n_2}^\pi(r)\mathbf{1}||_2+||Z^\pi||_{2}(||(\varphi^\pi(r)-\tilde{\varphi}^\pi(r))\mathbf{1}||_2+||\left(\tilde{\varphi}^\pi(r)-\varphi_{n_2}^\pi(r))\mathbf{1}||_2\right)
			\end{aligned}
		\end{equation*}	
		
		For $x\in\mathcal{X}$,
		\begin{equation*}
			\begin{aligned}|r(x)-\varphi_{n_2}^\pi(r)|=&\bigg|r(x)-\sum_{x'\in\mathcal{X}}\xi_{n_2}^\pi(x')r(x')\bigg|\\
				\leq&|r(x)|+\sum_{x'\in\mathcal{X}}\xi_{n_2}^\pi(x')|r(x')|\\
				\leq&|\mathcal{X}|+1,
			\end{aligned}
		\end{equation*}
		we have
		$||r-\varphi_{n_2}^\pi(r)\mathbf{1}||_{\infty}\leq|\mathcal{X}|+1.$
		Thus
		$||r-\varphi_{n_2}^\pi(r)\mathbf{1}||_2\leq\sqrt{|\mathcal{X}|}(|\mathcal{X}|+1)$.
		
		\eqref{invariant probabilities convergence} implies
		\begin{equation*}
			\mathbb{P}\left(||\tilde{\xi}^\pi-\xi^\pi||_2\leq\frac{\epsilon_1(n_1)}{1-\tau_1(P^\pi)}\right)\geq1-\delta_1(n_1).
		\end{equation*}
		Since $||r||_2\leq\sqrt{|\mathcal{X}|}$ and
		$|\tilde{\varphi}^\pi(r)-\varphi^\pi(r)|\leq||\tilde{\xi}^\pi-\xi^\pi||_2||r||_2$,
		we have
		$$\mathbb{P}\left(|\tilde{\varphi}^\pi(r)-\varphi^\pi(r)|\leq\frac{\epsilon_1(n_1)\sqrt{|\mathcal{X}|}}{1-\tau_1(P^\pi)}\right)\geq1-\delta_1(n_1),$$
		which further yields
		$$\mathbb{P}\left(||(\tilde{\varphi}^\pi(r)-\varphi^\pi(r))\mathbf{1}||_2\leq\frac{\epsilon_1(n_1)|\mathcal{X}|}{1-\tau_1(P^\pi)}\right)\geq1-\delta_1(n_1)$$
		due to the fact that
		$||(\tilde{\varphi}^\pi(r)-\varphi^\pi(r))\mathbf{1}||_2\leq\sqrt{|\mathcal{X}|}||(\tilde{\varphi}^\pi(r)-\varphi^\pi(r))\mathbf{1}||_\infty$.
		
		Moreover, \eqref{probabilistic convergence rate} implies
		\begin{equation*}
			\mathbb{P}\left(||\xi^\pi_{n_2}-\tilde{\xi}^\pi||_2\leq2M(x)\tau^{n}\right)\geq1-\delta_1(n_1),
		\end{equation*}
		Similarly, we have
		$$\mathbb{P}\left(||(\tilde{\varphi}^\pi(r)-\varphi_{n_2}^\pi(r))\mathbf{1}||_2\leq2M(x)\tau^{n}|\mathcal{X}|\right)\geq1-\delta_1(n_1).$$
		
		Therefore, we obtain
		\begin{equation*}
			\mathbb{P}\left(||\hat{r}^\pi-\hat{r}_{n_1,n_2}^\pi||_2\leq\varsigma\right)\geq1-6\delta_1(n_1)
		\end{equation*}
		where
		$\varsigma=\frac{\sqrt{|\mathcal{X}|}(|\mathcal{X}|+1)\widetilde{\varpi}}{c(c-\widetilde{\varpi})}+\frac{|\mathcal{X}|}{c}\left(\frac{\epsilon_1(n_1)}{1-\tau_1(P^\pi)}+2M(x)\tau^{n}\right)$. Finally,
		due to Lemma \ref{lem_ergodicity coefficient}, we have
		\begin{equation*}
			\mathbb{P}\left(||\hat{r}^\pi-\hat{r}_{n_1,n_2}^\pi||_2\leq\tilde{\varsigma}\right)\geq1-7\delta_1(n_1),
		\end{equation*}
		where
		$\tilde{\varsigma}=\frac{\sqrt{|\mathcal{X}|}(|\mathcal{X}|+1)\widetilde{\varpi}}{c(c-\widetilde{\varpi})}+\frac{|\mathcal{X}|}{c}\left(\widetilde{\varpi}-\epsilon_1(n_1)\sqrt{|\mathcal{X}|}\right)$
		and $\widetilde{\varpi}$ is defined in \eqref{def_tilde w}.
	\end{proof}
	\subsection{Proof of Theorem \ref{thm_N}}
	\begin{proof}
		Recall that the $x^{th}$ element of vector $P^\pi\hat{r}^\pi$
		is denoted by $P^\pi\hat{r}^\pi(x)$. Define the vectors
		$(\hat{r}^\pi)^2$ and $(P^\pi\hat{r}^\pi)^2$ with their
		$x^{th}$ elements being $(\hat{r}^\pi)^2(x)=\hat{r}^\pi(x)^2$
		and $(P^\pi\hat{r}^\pi)^2(x)=(P^\pi\hat{r}^\pi(x))^2$,
		respectively. Similarly for the vectors
		$(\hat{r}_{n_1,n_2}^\pi)^2$ and
		$(P_{n_1}^\pi\hat{r}_{n_1,n_2}^\pi)^2$.  Lemma
		\ref{lem_hat(r)} implies
		\begin{equation*}
			\begin{aligned}
				&||(\hat{r}^\pi)^2||_\infty\leq\frac{|\mathcal{X}|(|\mathcal{X}|+1)^2}{c^2},\\
				&||(P^\pi\hat{r}^\pi)^2||_\infty\leq\frac{|\mathcal{X}|^2(|\mathcal{X}|+1)^2}{c^2},
			\end{aligned}
		\end{equation*}
		and from Lemma \ref{lem_hat(r)_n} we get
		\begin{equation}\label{eq_hat(r)_n_norm}
			\begin{aligned}	
				&\mathbb{P}\left(||(\hat{r}_{n_1,n_2}^\pi)^2||_\infty\leq\frac{|\mathcal{X}|(|\mathcal{X}|+1)^2}{(c-\widetilde{\varpi})^2}\right)\geq1-4\delta_1(n_1),\\
				&\mathbb{P}\left(||(P^\pi_{n_1}\hat{r}_{n_1,n_2}^\pi)^2||_\infty\leq\frac{|\mathcal{X}|^2(|\mathcal{X}|+1)^2}{(c-\widetilde{\varpi})^2}\right)\geq1-4\delta_1(n_1).
			\end{aligned}
		\end{equation}
		First, by equations \eqref{asymptotic variance} and
		\eqref{estimated asymptotic variance}, we have
		\begin{equation*}
			\begin{aligned}	
				&|\sigma_{n_1,n_2}^2-\sigma^2|\\
				=&\bigg|\sum_{x\in\mathcal{X}}[{\hat{r}_{n_1,n_2}^\pi}(x)^2-(P_{n_1}^\pi\hat{r}_{n_1,n_2}^\pi(x))^2]\xi^\pi_{n_2}(x)-\sum_{x\in\mathcal{X}}[{\hat{r}^\pi}(x)^2-(P^\pi\hat{r}^\pi(x))^2]\xi^\pi(x)\bigg|\\
				\stackrel{(1)}{=}&\bigg|\sum_{x\in\mathcal{X}}[{\hat{r}_{n_1,n_2}^\pi}(x)^2-(P_{n_1}^\pi\hat{r}_{n_1,n_2}^\pi(x))^2]\xi^\pi_{n_2}(x)-\sum_{x\in\mathcal{X}}[{\hat{r}_{n_1,n_2}^\pi}(x)^2-(P_{n_1}^\pi\hat{r}_{n_1,n_2}^\pi(x))^2]\xi^\pi(x)\\
				&+\sum_{x\in\mathcal{X}}[{\hat{r}_{n_1,n_2}^\pi}(x)^2-(P_{n_1}^\pi\hat{r}_{n_1,n_2}^\pi(x))^2]\xi^\pi(x)-\sum_{x\in\mathcal{X}}[{\hat{r}^\pi}(x)^2-(P^\pi\hat{r}^\pi(x))^2]\xi^\pi(x)\bigg|\\
				&\stackrel{(2)}{\leq}a_n+b_n,
			\end{aligned}
		\end{equation*}
		where \begin{equation*}
			\begin{aligned}
				a_n=&\bigg|\sum_{x\in\mathcal{X}}[{\hat{r}_{n_1,n_2}^\pi}(x)^2-(P_{n_1}^\pi\hat{r}_{n_1,n_2}^\pi(x))^2](\xi^\pi_{n_2}(x)-\xi^\pi(x))\bigg|,\\
				b_n=&\bigg|\sum_{x\in\mathcal{X}}[{\hat{r}_{n_1,n_2}^\pi}(x)^2-(P_{n_1}^\pi\hat{r}_{n_1,n_2}^\pi(x))^2]\xi^\pi(x)-\sum_{x\in\mathcal{X}}[{\hat{r}^\pi}(x)^2-(P^\pi\hat{r}^\pi(x))^2]\xi^\pi(x)\bigg|.
			\end{aligned}
		\end{equation*}
		Equality $(1)$ holds due to subtracting and then adding a
		same term
		$\sum_{x\in\mathcal{X}}[{\hat{r}_{n_1,n_2}^\pi}(x)^2-(P_{n_1}^\pi\hat{r}_{n_1,n_2}^\pi(x))^2]\xi^\pi(x)$. Inequality
		$(2)$ holds by triangle inequality. In the following, we
		bound the terms $a_n$ and $b_n$ individually.  Term $a_n$ can
		be further expressed as
		\begin{equation*}
			\begin{aligned}
				a_n\stackrel{(1)}{\leq}&\sum_{x\in\mathcal{X}}|[{\hat{r}_{n_1,n_2}^\pi}(x)^2-(P_{n_1}^\pi\hat{r}_{n_1,n_2}^\pi(x))^2](\xi^\pi_{n_2}(x)-\xi^\pi(x))|\\
				\stackrel{(2)}{\leq}&||({\hat{r}_{n_1,n_2}^\pi})^2-(P_{n_1}^\pi\hat{r}_{n_1,n_2}^\pi)^2||_\infty||\xi^\pi_{n_2}-\xi^\pi(\cdot)||_1\\
				\stackrel{(3)}{\leq}&\left(||({\hat{r}_{n_1,n_2}^\pi})^2||_\infty+||(P_{n_1}^\pi\hat{r}_{n_1,n_2}^\pi)^2||_\infty\right)\left(||\xi^\pi_{n_2}-\tilde{\xi}^\pi||_1+||\tilde{\xi}^\pi-\xi^\pi||_1\right).
			\end{aligned}
		\end{equation*}
		
		Inequalities $(1)$ and $(3)$ hold due to triangle property,
		while $(2)$ is true because of H\"{o}lder's inequality. From
		\eqref{probabilistic convergence rate} and \eqref{invariant
			probabilities convergence}, we have
		\begin{equation}\label{an_1}
			\mathbb{P}\left(||\xi^\pi_{n_2}-\tilde{\xi}^\pi||_1+||\tilde{\xi}^\pi-\xi^\pi||_1\leq\iota\right)\geq1-2\delta_1(n_1).
		\end{equation}
		where
		$\iota=2\sqrt{|\mathcal{X}|}M(x)\tau^{n_2}+\frac{\epsilon_1(n_1)}{1-\tau_1(P^\pi)}$. In
		addition, under Lemma \ref{lem_ergodicity coefficient} we
		obtain
		\begin{equation*}
			\mathbb{P}\left(||\xi^\pi_{n_2}-\tilde{\xi}^\pi||_1+||\tilde{\xi}^\pi-\xi^\pi||_1\leq\tilde{\iota}\right)\geq1-3\delta_1(n_1),
		\end{equation*}		
		where
		$\tilde{\iota}=2\sqrt{|\mathcal{X}|}M(x)\tau^{n_2}+\frac{\epsilon_1(n_1)}{1-\tau_1(P_{n_1}^\pi)-\epsilon_1(n_1)}$. Combining
		\eqref{eq_hat(r)_n_norm} and \eqref{an_1}, by the union bound
		of probability and definition of $\widetilde{\varpi}$, we
		have
		\begin{equation*}
			\begin{aligned}
				\mathbb{P}&\left(a_n\leq\frac{|\mathcal{X}|(|\mathcal{X}|+1)^3(\widetilde{\varpi}-\sqrt{|\mathcal{X}|}\epsilon_1(n_1))}{(c-\widetilde{\varpi})^2}\right)\leq1-11\delta_1(n_1).
			\end{aligned}
		\end{equation*}
		Next, we bound the term $b_n$.
		\begin{equation*}
			\begin{aligned}
				b_n\stackrel{(1)}{\leq}&\sum_{x\in\mathcal{X}}|\{[{\hat{r}_{n_1,n_2}^\pi}(x)^2-(P_{n_1}^\pi\hat{r}_{n_1,n_2}^\pi(x))^2]-[{\hat{r}^\pi}(x)^2-(P^\pi\hat{r}^\pi(x))^2]\}\xi^\pi(x)|\\
				\stackrel{(2)}{\leq}&||(\hat{r}^\pi_{n_1,n_2})^2-(\hat{r}^\pi)^2||_\infty||\xi^\pi||_1+||(P^\pi_{n_1}\hat{r}^\pi_{n_1,n_2})^2-(P^\pi\hat{r})^2||_\infty||\xi^\pi||_1\\
				\stackrel{(3)}{\leq}&|\mathcal{X}|(||(\hat{r}^\pi_{n_1,n_2})^2-(\hat{r}^\pi)^2||_\infty+||(P^\pi_{n_1}\hat{r}^\pi_{n_1,n_2})^2-(P^\pi\hat{r})^2||_\infty).
			\end{aligned}
		\end{equation*}
		Inequality $(1)$ holds due to triangle property, and $(2)$ is
		due to H\"{o}lder's inequality, while $(3)$ is true since
		$\xi^\pi(x)\in[0,1]$.  Note that
		\begin{equation*}
			\begin{aligned}
				&|{\hat{r}_{n_1,n_2}^\pi}(x)^2-{\hat{r}^\pi}(x)^2|\\
				=&|{\hat{r}_{n_1,n_2}^\pi}(x)-{\hat{r}^\pi}(x)||{\hat{r}_{n_1,n_2}^\pi}(x)+{\hat{r}^\pi}(x)|\\
				\leq&|{\hat{r}_{n_1,n_2}^\pi}(x)-{\hat{r}^\pi}(x)|(|{\hat{r}_{n_1,n_2}^\pi}(x)|+|{\hat{r}^\pi}(x)|).
			\end{aligned}
		\end{equation*}	
		Furthermore, from Lemma \ref{lem_hat(r)}, \ref{lem_hat(r)_n}
		and \ref{lem_solution difference}, we have
		\begin{equation}\label{bn-1}
			\begin{aligned}
				\mathbb{P}\left(||(\hat{r}^\pi_{n_1,n_2})^2-(\hat{r}^\pi)^2||_\infty\leq\vartheta\right)\geq1-11\delta_1(n_1).
			\end{aligned}
		\end{equation}
		where \begin{equation*}
			\begin{aligned}
				\vartheta=&\frac{|\mathcal{X}|(|\mathcal{X}|+1)(2c-\widetilde{\varpi})}{c^2(c-\widetilde{\varpi})}\left[\frac{(|\mathcal{X}|+1)\widetilde{\varpi}}{c-\widetilde{\varpi}}+\sqrt{|\mathcal{X}|}(\widetilde{\varpi}-\epsilon_1(n_1)\sqrt{|\mathcal{X}|})\right].
			\end{aligned}
		\end{equation*}	
		In addition, we have
		\begin{equation*}
			\begin{aligned}
				&|(P_{n_1}^\pi\hat{r}_{n_1,n_2}^\pi(x))^2-(P^\pi\hat{r}^\pi(x))^2|
				\\=&|P_{n_1}^\pi\hat{r}_{n_1,n_2}^\pi(x)-P^\pi\hat{r}^\pi(x)|\times|P_{n_1}^\pi\hat{r}_{n_1,n_2}^\pi(x)+P^\pi\hat{r}^\pi(x)|\\
				\stackrel{(1)}{\leq}&(|P_{n_1}^\pi\hat{r}_{n_1,n_2}^\pi(x)-P^\pi\hat{r}_{n_1,n_2}^\pi(x)|+|P^\pi\hat{r}_{n_1,n_2}^\pi(x)-P^\pi\hat{r}^\pi(x)|)(|P_{n_1}^\pi\hat{r}_{n_1,n_2}^\pi(x)|+|P^\pi\hat{r}^\pi(x)|)\\
				=&\left(\bigg|\sum_{x'\in\mathcal{X}}(P^\pi_{n_1}(x,x')-P^\pi(x,x'))\hat{r}^\pi_{n_1,n_2}(x')|+|\sum_{x'\in\mathcal{X}}P^\pi(x,x')(\hat{r}^\pi_{n_1,n_2}(x')-\hat{r}^\pi(x'))\bigg|\right)\\
				&\times(|P_{n_1}^\pi\hat{r}_{n_1,n_2}^\pi(x)|+|P^\pi\hat{r}^\pi(x)|)\\
				\stackrel{(2)}{\leq}&\left(\sum_{x'\in\mathcal{X}}|(P^\pi_{n_1}(x,x')-P^\pi(x,x'))\hat{r}^\pi_{n_1,n_2}(x')|+\sum_{x'\in\mathcal{X}}|P^\pi(x,x')(\hat{r}^\pi_{n_1,n_2}(x')-\hat{r}^\pi(x'))|\right)\\
				&\times(|P_{n_1}^\pi\hat{r}_{n_1,n_2}^\pi(x)|+|P^\pi\hat{r}^\pi(x)|)\\
				\stackrel{(3)}{\leq}&b_1\times(|P_{n_1}^\pi\hat{r}_{n_1,n_2}^\pi(x)|+|P^\pi\hat{r}^\pi(x)|).\\
			\end{aligned}
		\end{equation*}
		where
		\begin{equation*}
			\begin{aligned}
				b_1=||P^\pi_{n_1}(x,\cdot)-P^\pi(x,\cdot)||_1||\hat{r}^\pi_{n_1,n_2}||_\infty+||P^\pi(x,\cdot)||_1||\hat{r}^\pi_{n_1,n_2}-\hat{r}^\pi||_\infty.
			\end{aligned}
		\end{equation*}
		Inequalities $(1)$ and $(2)$ are due to triangle property,
		and $(3)$ is because of H\"{o}lder's inequality.
		
		Since
		$||P^\pi_{n_1}(x,\cdot)-P^\pi(x,\cdot)||_1\leq||P^\pi_{n_1}-P^\pi||_\infty$,
		Assumption \ref{asu_estimated probability matrix} implies
		\begin{equation*}
			\mathbb{P}\left(||P^\pi_{n_1}(x,\cdot)-P^\pi(x,\cdot)||_1\leq\epsilon_1(n_1)\right)\geq 1-\delta_1(n_1).
		\end{equation*}
		Since $||P^\pi(x,\cdot)||_1\leq|\mathcal{X}|$, from Lemma
		\ref{lem_hat(r)_n} and \ref{lem_solution difference}, we have
		\begin{equation*}
			\begin{aligned}
				\mathbb{P}\left(b_1\leq\overline{b}_1\right)\geq
				1-12\delta_1(n_1).
			\end{aligned}
		\end{equation*}
		where
		\begin{equation*}
			\begin{aligned}
				\overline{b}_1=&\frac{\epsilon_1(n_1)\sqrt{|\mathcal{X}|}(|\mathcal{X}|+1)}{c-\widetilde{\varpi}}+\frac{|\mathcal{X}|^{\frac{3}{2}}}{c}\left(\frac{\widetilde{\varpi}(|\mathcal{X}|+1)}{c-\widetilde{\varpi}}+\sqrt{|\mathcal{X}|}\left(\widetilde{\varpi}-\epsilon_1(n_1)\sqrt{|\mathcal{X}|}\right)\right).
			\end{aligned}
		\end{equation*}
		Moreover, from Lemma \ref{lem_hat(r)} and \ref{lem_hat(r)_n},
		we obtain
		\begin{equation*}
			\mathbb{P}\left(|P_{n_1}^\pi\hat{r}_{n_1,n_2}^\pi(x)|+|P^\pi\hat{r}^\pi(x)|\leq\eta\right)\geq 1-4\delta_1(n_1)
		\end{equation*}
		where
		$\eta=\dfrac{|\mathcal{X}|(|\mathcal{X}|+1)(2c-\widetilde{\varpi})}{c(c-\widetilde{\varpi})}$.
		Note $\widetilde{\varpi}$ is dependent on the initial state
		$x$ and denote
		\begin{equation*}
			\begin{aligned}
				\alpha_1(x)=&\frac{\epsilon_1(n_1)\sqrt{|\mathcal{X}|}(|\mathcal{X}|+1)}{c-\widetilde{\varpi}}+\frac{|\mathcal{X}|^{\frac{3}{2}}}{c}\left(\frac{\widetilde{\varpi}(|\mathcal{X}|+1)}{c-\widetilde{\varpi}}+\sqrt{|\mathcal{X}|}\left(\widetilde{\varpi}-\epsilon_1(n_1)\sqrt{|\mathcal{X}|}\right)\right),\\
				\alpha_2(x)=&\dfrac{|\mathcal{X}|(|\mathcal{X}|+1)(2c-\widetilde{\varpi})}{c(c-\widetilde{\varpi})},
			\end{aligned}
		\end{equation*}
		we conclude
		\begin{equation}\label{bn-2}
			\begin{aligned}
				&\mathbb{P}\left(||(P_{n_1}^\pi\hat{r}_{n_1,n_2}^\pi)^2-(P^\pi\hat{r}^\pi)^2||_\infty\leq\max_{x\in\mathcal{X}}\alpha_1(x)\alpha_2(x)\right)\geq1-16\delta_1(n_1).
			\end{aligned}
		\end{equation}
		
		At last, by the union bound of probability, \eqref{bn-1} and
		\eqref{bn-2} yield
		\begin{equation*}
			\begin{aligned}
				\mathbb{P}\left(b_n\leq\overline{b}_n\right)\geq1-27\delta_1(n_1),
			\end{aligned}
		\end{equation*}
		where
		\begin{equation*}
			\begin{aligned}
				\overline{b}_n&=\frac{\alpha_2(x)|\mathcal{X}|}{c}\left[\frac{(|\mathcal{X}|+1)\widetilde{\varpi}}{c-\widetilde{\varpi}}+\sqrt{|\mathcal{X}|}(\widetilde{\varpi}-\epsilon_1(n_1)\sqrt{|\mathcal{X}|})\right]+\max_{x\in\mathcal{X}}\alpha_1(x)\alpha_2(x)|\mathcal{X}|.
			\end{aligned}
		\end{equation*}
		Recall that for two random variables $a_n$ and $b_n$, and
		$\epsilon\in\mathbb{R}$ with $\lambda\in(0,1)$, by union
		bound and monotonicity property of probability, we have
		$$\mathbb{P}(a_n+b_n\leq \epsilon)\geq\mathbb{P}(a_n\leq\lambda \epsilon)+\mathbb{P}(b_n\leq(1-\lambda) \epsilon).$$
		Therefore, if for $a_n$
		$$\frac{|\mathcal{X}|(|\mathcal{X}|+1)^3(\widetilde{\varpi}-\sqrt{|\mathcal{X}|}\epsilon_1(n_1))}{(c-\widetilde{\varpi})^2}\leq\lambda\epsilon,$$
		and for $b_n$
		\begin{equation*}
			\begin{aligned}
				&\frac{\alpha_2(x)|\mathcal{X}|}{c}\left[\frac{(|\mathcal{X}|+1)\widetilde{\varpi}}{c-\widetilde{\varpi}}+\sqrt{|\mathcal{X}|}(\widetilde{\varpi}-\epsilon_1(n_1)\sqrt{|\mathcal{X}|})\right]+\max_{x\in\mathcal{X}}\alpha_1(x)\alpha_2(x)|\mathcal{X}|\leq(1-\lambda)\epsilon.
			\end{aligned}
		\end{equation*}
		
		We can conclude that
		$$\mathbb{P}(a_n+b_n\leq \epsilon)\geq1-38\delta_1(n_1),$$
		which further implies
		$$\mathbb{P}(|\sigma_{n_1,n_2}^2-\sigma^2|\leq\epsilon)\geq1-38\delta_1(n_1)$$
		by the monotonicity property and union bound of probability.
	\end{proof}
	\subsection{Proof of Corollary \ref{col_2}}
	\begin{proof}
		From the definitions of $\tilde{\rho}_\lambda(S^\pi_T)$ and
		$\rho_\lambda(S^\pi_T)$, we get
		\begin{equation*}
			\begin{aligned}
				&|\tilde{\rho}_\lambda(S^\pi_T)-\rho_\lambda(S^\pi_T)|\\
				=&|\varphi(r)-\varphi_{n_2}^\pi(r)+\lambda\sigma_{n_1,n_2}^2-\lambda\sigma^2|\\
				\stackrel{(1)}{\leq}&\bigg|\sum_{x\in\mathcal{X}}[\xi^\pi(x)-\xi^\pi_{n_2}(x)]r(x)\bigg|+\lambda|\tilde{\sigma}^2-\sigma^2|\\
				\stackrel{(2)}{\leq}&\sum_{x\in\mathcal{X}}[|\xi^\pi(x)-\tilde{\xi}^\pi(x)|+|\tilde{\xi}^\pi(x)-\xi^\pi_{n_2}(x)|]|r(x)|+\lambda|\sigma_{n_1,n_2}^2-\sigma^2|\\
				\stackrel{(3)}{\leq}&||\xi^\pi-\tilde{\xi}^\pi||_1+||\tilde{\xi}^\pi-\xi^\pi_{n_2}||_1+\lambda|\sigma_{n_1,n_2}^2-\sigma^2|.
			\end{aligned}
		\end{equation*}
		Here inequalities $(1)$ and $(2)$ are due to triangle
		property. Inequality $(3)$ is because H\"{o}lder's inequality
		and $||r||_\infty\leq1$.  Moreover, since $P^\pi_{n_1}$ is the
		{\em $\epsilon$-perturbed transition kernel} under Assumption
		\ref{asu_convergence rate}, by \eqref{probabilistic
			convergence rate} and \eqref{estimated stationary
			distribution}, we have
		\begin{equation}\label{col_2_1}
			\mathbb{P}\left(||\tilde{\xi}^\pi-\xi^\pi_{n_2}||_1\leq 2M(x)\tau^{n_2}\right)\geq1-\delta_1(n_1).
		\end{equation}		
		In addition, from \eqref{invariant probabilities convergence}
		and Lemma \ref{lem_ergodicity coefficient}, we get
		\begin{equation}\label{col_2_2}
			\mathbb{P}\left(||\tilde{\xi}^\pi-\xi^\pi||_1\leq\frac{\epsilon_1(n_1)}{1-\tau_1(P_{n_1}^\pi)-\epsilon_1(n_1)}\right)\geq1-2\delta_1(n_1).
		\end{equation}	
		Applying the union bound of probability to \eqref{col_2_1},
		\eqref{col_2_2} and Theorem \ref{thm_N}, we have
		\begin{equation*}
			\mathbb{P}\left(|\tilde{\rho}_\lambda(S^\pi_T)-\rho_\lambda(S^\pi_T)|\leq\kappa\right)\geq1-41\delta_1(n_1).
		\end{equation*}
		where
		$\kappa=2M(x)\tau^{n_2}+\frac{\epsilon_1(n_1)}{1-\tau_1(P_{n_1}^\pi)-\epsilon_1(n_1)}+\lambda\epsilon$.
	\end{proof}
	\subsection{Proof of Theorem \ref{thm_improvement}}
	\begin{proof}
		Note that, for fixed $\theta$, we are able to evaluate the
		policy $\rho(S^\theta_T)$. Moreover, \eqref{update} can be
		considered as a discretization of the differential equation
		\eqref{differential equation}.
		
		For SPSA update, $\mathcal{Z}$ is an asymptotically stable
		attractor for the differential equation \eqref{differential
			equation}, with $\rho(S_T^\theta)$ itself serving as a
		strict Lyapunov function. This can be inferred as follows
		$$\frac{d\rho(S_T^\theta)}{dt}=\nabla_{\theta}\rho(S_T^\theta)\dot{\theta}=\nabla_{\theta}\rho(S_T^\theta)\check{\Gamma}(-\nabla_{\theta}\rho(S_T^\theta))<0.$$ The claim for SPSA update now follows from \citep[Theorem 5.3.1]{kushner2012stochastic}.
		
		The convergence analysis for SF update follows in a similar
		manner as SPSA update.
	\end{proof}
	
	\bibliographystyle{plainnat} 
	\bibliography{reference}
		
\end{document}